\documentclass{amsart}

\usepackage{amsmath, amssymb,amsthm}%, 
\usepackage{inputenc, graphicx, verbatim,color, enumerate}

\usepackage{float}

\usepackage{tikz,pgfplots}
\usepackage{url}
\usepackage{xcolor}

\newcommand{\hand}{\hbox{ and }}

\newcommand{\F}{\mathbb{F}}

\newcommand{\Z}{\mathbb{Z}}
\renewcommand{\P}{\mathbb{P}}

\newcommand{\UU}{\mathcal{U}}
\newcommand{\FF}{\mathcal{F}}
\newcommand{\ZZ}{\mathcal{Z}}
\newcommand{\CC}{\mathcal{C}}
\newcommand{\cyc}{\mathrm{cyc}}
\newcommand{\cl}{\mathrm{cl}}

\newcommand{\supp}{\mathrm{supp}}
\newcommand{\wt}{\mathrm{wt}}
\newcommand{\ie}{\emph{i.e.}}

\newcommand{\nd}{\; \mathrm{and} \; }
\newcommand{\st}{\; \mathrm{such\ that} \;}

\newcommand{\rmatroid}{M = (E,\rho)}

\newtheorem{theorem}{Theorem}
\newtheorem{proposition}{Proposition}
\newtheorem{lemma}{Lemma}
\newtheorem{corollary}{Corollary}
\theoremstyle{definition}\newtheorem{definition}{Definition}
\theoremstyle{definition}\newtheorem{example}{Example}
\theoremstyle{definition}\newtheorem{notation}{Notation}

\begin{document}

%\begin{frontmatter}

\title{Cyclic Flats of Binary Matroids}

\author[aalto]{Ragnar Freij-Hollanti}
\author[aalto]{Matthias Grezet}
\author[aalto]{Camilla Hollanti} 
\author[malardalen]{Thomas Westerb\"ack}
\address[Freij-Hollanti, Grezet and Hollanti]{Department of Mathematics and Systems Analysis, Aalto University, FI-00076 Aalto, Finland}
\address[Westerb\"ack]{Division of Applied Mathematics, UKK, M\"{a}lardalen University, H\"ogskoleplan 1, Box 883, 721 23 V\"aster\r{a}s, Sweden}

\maketitle

\begin{abstract}
In this paper, first steps are taken towards characterising lattices of cyclic flats $\ZZ(M)$ that belong to matroids $M$ that can be represented over a prescribed finite field $\F_q$. Two natural maps from $\ZZ(M)$ to the lattice of cyclic flats of a minor of $M$ are given. Binary matroids are characterised via their lattice of cyclic flats. It is shown that the lattice of cyclic flats of a simple binary matroid without isthmuses is atomic.

\end{abstract}

%\begin{keyword}
%Cyclic flats \sep Binary matroids \sep Atomic lattices

 %\MSC 05B35 
%\end{keyword}

%\end{frontmatter}

%\tableofcontents

%Input "intro"
%Contains : Introduction + Preliminaries 

%Matroids were first introduced by Whitney in 1935, to capture and generalise the notion of linear dependence in purely combinatorial terms. %~\cite{whitney35}. 
%Indeed, the combinatorial setting is general enough to also capture many other notions of dependence occurring in mathematics, such as cycles or incidences in a graph, non-transversality of algebraic varieties, or algebraic dependence of field extensions. Of special interest for linear LRCs is the connection between linear algebra and matroids. /Users/cash/Dropbox/Matroid research/Matthias/BinaryLCF journal/BinaryLCF journal/input1_intro.tex
%For an introduction to matroid theory, see for example \cite{oxley11}. 

\section{Introduction}

In traditional matroid theory, one of the most crucial objects is that of a lattice of flats. This is a geometric lattice, \ie,  it is atomic and semimodular, and in fact every geometric lattice is the lattice of flats $\FF(M)$ of a simple matroid $M=(E,\rho)$~\cite{birkhoff1995}. This correspondence between lattices and matroids behaves reasonably well with respect to their respective notions of duality, namely, the dual lattice of $\FF(M)$ is isomorphic to the lattice of cyclic sets $\UU(M^*)$, whose elements are unions of circuits in the dual $M^*$.

Thus, the Boolean lattice $2^E$ has two subposets $\FF(M)$ and $\UU(M)$, both of which are lattices, each of which determine the matroid $M$ uniquely. This has inspired many authors to look at their intersection $\ZZ(M)=\FF(M)\cap\UU(M)$~\cite{sims80,bonin08,prideaux16, eberhardt14}. It was shown independently in~\cite{sims80,bonin08} that $\ZZ(M)$ together with the restriction of the rank function $\rho$ to $\ZZ(M)$ is enough to determine $M$. Moreover, $\ZZ(M)$ is a lattice, although its lattice structure is neither induced by $2^E$, $\FF(M)$, or $\UU(M)$~\cite{sims80,bonin08}.

As opposed to $\FF(M)$ and $\UU(M)$, the lattice of cyclic flats has no additional structure apart from being a lattice. Indeed, it is shown in~\cite{bonin08} that every finite lattice is isomorphic to the lattice of cyclic flats of some finite matroid. Yet, there are many advantages in describing a matroid in terms of its lattice of cyclic flats. Firstly, the cyclic flats description is rather concise for many naturally occurring matroids. Secondly, it was shown in~\cite{westerback15} that many central invariants in coding theory can be naturally described in terms of the lattice of cyclic flats of the associated matroid. Especially, this was shown to be the case for invariants related to applications to distributed data storage~\cite{westerback15}. It was earlier shown in~\cite{eberhardt14} that the Tutte polynomial can be computed efficiently for matroids whose lattice of cyclic flats has bounded height. Yet another reason to take interest in the lattices of cyclic flats is that some natural classes of matroids can be defined in terms of the structure of $\ZZ(M)$. For instance, a matroid $M$ is nested if and only if $\ZZ(M)$ is a chain, The nested matroids form the first known example of a minor-closed class of matroids that is well-quasi-ordered under the minor relation, but has infinitely many forbidden minors.

In this work, we are taking first steps towards characterising lattices of cyclic flats that belong to matroids that can be represented over a prescribed finite field $\F_q$. Our approach is to study the minor relation from the point of view of cyclic flats. In particular, in Theorem~\ref{thm:cf_formulas} we construct two natural maps from $\ZZ(M)$ to the lattice of cyclic flats of a minor of $M$.

We take inspiration from Rota's conjecture~\cite{rota71}, and its recently announced proof~\cite{whittle14}, that representability over a prescribed finite field $\F_q$ is equivalent to avoiding a finite set $L(\F_q)$ of minors. However, in this initial work we only actually use the rather weak result that if $n>q+1$, and the uniform matroid $U_n^2$ is a minor of $M$, then $M$ is not representable over $\F_q$~\cite{oxley11}. Thus, in Theorem~\ref{thm:u2n}, we compute the largest $n$ for which $U_n^2$ is a minor of $M$, from $\ZZ(M)$. This is done via studying a certain antichain of flats associated to every cyclic flat of rank $\rho(1_{\ZZ(M)})-2$. By duality, of course, this can also be used to find the largest $n$ for which $U_n^{n-2}$ is a minor of $M$.

The representability over $\F_{2}$, or other small fields of characteristic $2$, is particularly interesting from a data storage point of view. For instance, small fields allow an efficient implementation of locally repairable codes \cite{papailiopoulos12}. Constructions of optimal locally repairable codes over $\F_{2}$ were also derived in \cite{silberstein18,huang15}. Therefore, it motivates a deeper understanding of the dependency structures of binary matroids. 

Since binary matroids are exactly characterised by not having $U_4^2$ as a minor, we thus get two equivalent necessary and sufficient conditions for a matroid to be representable over $\F_2$ in Corollary~\ref{cor:binary}. In the second half of the paper, we focus exclusively on binary matroids. In Section 7 and 8, we study sublattices of $\ZZ(M)$ of height 2 and 3 respectively, when $M$ is binary. We also prove, in Theorem~\ref{thm:lcf_atomic}, that the lattice of cyclic flats of a simple matroid with no isthmuses is atomic.

Understanding the sublattices of small height helps us describe constraints on $\ZZ(M)$ recursively in Section 9, and these recursive constraints are enough to reprove the Griesmer bound for binary codes~\cite{huffman2010}. On our way to proving the Griesmer bound, we define the class of {\em blunt} cyclic flats of a binary matroid. These play a special role in our analysis and seem relevant also in a much broader context, although it is not clear how to generalise the definition to non-binary matroids.

Part of this work has previously been presented at the $5^{\rm th}$ International Castle Meeting on Coding Theory and Applications~\cite{grezet17} and at the International Zurich Seminar
on Information and Communication \cite{grezet18}.

\section{Preliminaries}

Matroids have many equivalent definitions in the literature. Here, we choose to present matroids via their rank functions. %Much of the contents in this section can be found in more detail in \cite{freij-hollanti17}.

\begin{definition}\label{def:matroid_rank}
A \emph{(finite) matroid} $\rmatroid$ is a finite set $E$ together with a \emph{rank function} $\rho:2^E \rightarrow \Z$ such that for all subsets $X,Y \subseteq E$
\begin{itemize}
\item[(R1)] $0 \leq \rho(X) \leq |X|$,
\item[(R2)] If $X \subseteq Y$ then $\rho(X) \leq \rho(Y)$,
\item[(R3)] $\rho(X) + \rho(Y) \geq \rho(X \cup Y) + \rho(X \cap Y)$. 
\end{itemize}
\end{definition}

When $\rmatroid$ is a matroid, we also define the \emph{nullity function} $\eta:2^E\to\Z$ by $\eta(X)=|X|-\rho(X)$.

%A subset $X \subseteq E$ is called \emph{independent} if $\rho(X) = |X|$. If $X$ is independent and $\rho(X) = \rho(E)$, then $X$ is called a \emph{basis}. Strongly related to the rank function is the \emph{nullity function} $\eta:2^E \rightarrow \mathbb{Z}$, defined by $\eta(X) = |X| - \rho(X)$ for $X \subseteq E$.

Any matrix $G$ over a field $\F$ generates a matroid $M_G=(E,\rho)$, where $E$ is the set of columns of $G$, and $\rho(X)$ is the rank of $G(X)$ over $\F$, where $G(X)$ denotes the submatrix of $G$ formed by the columns indexed by $X$. Clearly, $\rho$ only depends on the row space of $G$, so row-equivalent matrices generate the same matroid.  

%Thus, there is a straightforward connection between linear codes and matroids. Let $C$ be a linear code over a field $\F_q$. %generated by $G$, meaning that $C$ equals the rowspace of $G$. 
%Then any two different generator matrices of $C$ will have the same row space by definition, so they will generate the same matroid. Therefore, without any inconsistency, we can denote the matroid associated to these generator matrices by $M_C = (E,\rho_C)$. The rank function $\rho_C$ can be defined directly from the code without referring to a generator matrix, via $\rho_C(X) = \mathrm{dim}(C|X)$ for $X \subseteq E$.

Two matroids $M_1 = (E_1,\rho_1)$ and $M_2 = (E_2,\rho_2)$ are \emph{isomorphic} if there exists a bijection $\psi: E_1 \rightarrow E_2$ such that $\rho_2(\psi(X)) = \rho_1(X)$ for all subsets $X \subseteq E_1$.

\begin{definition} A matroid that is isomorphic to $M_G$ for some matrix $G$ over $\F$ is said to be \emph{representable} over $\F$. %We also say that such a matroid is $\F$-representable. 
A \emph{binary} matroid is a matroid that is representable over $\F_2$.\end{definition}

\begin{definition}
The \emph{uniform matroid} $U_n^k=([n],\rho)$ is a matroid with a ground set $[n]=\{1,2,\ldots , n\}$ and a rank function $\rho(X)=\min\{|X|, k\}$ for $X\subseteq [n]$.  
\end{definition}

%Recall that a linear code $C\subseteq \F^n$ of length $n$ and dimension $k$ is said to be maximum distance separable (MDS) if its minimum Hamming distance satisfies $d=n-k+1$. Since such a code is characterised by having dimension $k$ when projected to a coordinate set of size $\geq k$, we get the following standard result. %The following straightforward observation gives a characterisation of maximum distance separable (MDS) codes and also shows that uniform matroids constitute a subclass of representable matroids.

%\begin{proposition}
%A linear code $C$ is an $[n,k,n-k+1]$-MDS code if and only if $M_C$ is the uniform matroid $U_n^k$.
%\end{proposition}

Motivated by coding theory and the relation between linear codes and matroids, we define the minimum distance of a matroid to be the following. 

\begin{definition}
Let $\rmatroid$ be a matroid. The minimum distance of $M$ is 
\[
d=\min\{ |X| : X \subseteq E, \rho(E-X)<\rho(E) \}.
\]
A matroid with $|E|=n$, $\rho(E)=k$, and minimum distance $d$ is referred to as an $(n,k,d)$-matroid. 
\end{definition}

Therefore, if the matroid $M_{\CC}$ comes from a linear code $\CC$, then the minimum distance of $M_{\CC}$ coincides with the minimum Hamming distance.

%There are several elementary operations that are useful for explicit constructions of matroids, as well as for analysing their structure. The operations that we will need for this paper are dualisation, contraction,  and deletion.

\begin{definition}
Let $\rmatroid$ be a matroid and $X, Y \subseteq E$, and denote by $\bar{X}=E-X$ for any $X\subseteq E$. Then
\begin{enumerate}[(i)]
\item The \emph{restriction} of $M$ to $Y$ is the matroid $M|Y = (Y,\rho_{|Y})$, where $\rho_{|Y}(A) = \rho(A)$ for $A \subseteq Y$. The restriction operation to $Y$ is also referred to as \emph{deletion} of the set $E-Y$. 
\item The \emph{contraction} of $M$ by $X$ is the matroid $M/X = (\bar{X},\rho_{/X})$, where $\rho_{/X}(A) = \rho(A \cup X) - \rho(X)$ for $A \subseteq \bar{X}$.
\item For $X \subseteq Y$, a \emph{minor} of $M$ is a matroid isomorphic to $M|Y/X = (Y-X,\rho_{|Y/X})$, obtained from $M$ by restriction to a set $Y\subseteq E$ and contraction by $X\subseteq Y$. Observe that this does not depend on the order in which the restriction and contraction are performed.
\item The \emph{dual} of $M$ is the matroid $M^* = (E, \rho^*)$, where 
\[
\rho^*(A) = |A| + \rho(\bar{A}) - \rho(E)=\eta(E)-\eta(\bar{A}) \text{ for } A \subseteq E.
\]
\end{enumerate}
\end{definition}

It is easy to see that representability over $\F_{q}$ is preserved under minors and duals. Given the structure of uniform matroids and the definition of a minor, the minors of uniform matroids are very easily described: 
\begin{lemma}\label{lm:subunif}
Let $U_n^k=([n], \rho)$ be a uniform matroid, and let $X\subseteq Y\subseteq E$. Then the minor $U_n^k|Y/X$ is isomorphic to $U_{n'}^{k'}$, where $k'=\max\{0,k-|X|\}$ and $n'=|Y|-|X|$. In particular, $M$ is a minor of $U_n^k$ if and only if $M\cong U_{n'}^{k'}$, for some $0\leq k'\leq k$ and $0\leq n'-k'\leq n-k$.
\end{lemma}
\begin{comment}
\begin{proof}
The first statement is an easy calculation from the definition of the rank function of the minor. For the second statement, notice that the possible values of $k'$ for a minor are precisely $0\leq k'\leq k$, and that and $n'-k'=\min\{|Y|-|X|, |Y|-k\}$, which can take any value between $0$ and $n-k$.
\qed\end{proof}
\end{comment}

In general there is no simple criterion to determine if a matroid is representable~\cite{vamos78,mayhew18}. However, there is a simple criterion for when a matroid is binary. 

\begin{theorem}[\cite{tutte58}]\label{tutte}
Let $\rmatroid$ be a matroid. The following two conditions are equivalent.
\begin{enumerate}
\item $M$ is representable over $\F_2$.
\item There are no sets $X\subseteq Y\subseteq E$ such that $M|Y/X$ is isomorphic to the uniform matroid $U_4^2$.
\end{enumerate}
\end{theorem}

If $M$ is representable over a fixed finie field $\F_{q}$, then so are all its minors. The class of matroids representable over $\F_{q}$ is therefore closed under minors. The following result, which extend the previous theorem, was first conjectured by Gian-Carlo Rota in 1970 \cite{rota71}. A proof of this conjecture was announced by Geelen, Gerards, and Whittle in 2014, but the details of the proof remain to be written up~\cite{whittle14}. 

\begin{theorem}[\hspace{1sp}\cite{whittle14}]
\label{thm:rota}
For any finite field $\F_{q}$, there is a finite set $L(\F_{q})$ of matroids such that any matroid $M$ is representable over $\F_{q}$ if and only if it contains no element from $L(\F_{q})$ as a minor.
\end{theorem}

Since the 1970's, it has been known that a matroid is representable over $\F_3$ if and only if it avoids the uniform matroids $U_5^2$, $U_5^3$, the Fano plane $\P^2(\F_2)$, and its dual $\P^2(\F_2)^*$ as minors. The list $L(\F_4)$ was given explicitly in 2000 and contains seven elements. For larger fields, the explicit list is not known, and there is little hope to even find useful bounds on its size. Assuming the MDS conjecture~\cite{Segre55}, a matroid $M$ that is linearly representable over $\F_q$ must avoid $U_{q+2}^k$ as a minor, for $k=2$, $4\leq k\leq q-2$, and $k=q$. If $q$ is odd, $M$ must also avoid $U_{q+2}^3$ and $U_{q+2}^{q-1}$ minors. The MDS conjecture is widely believed to be true, and is proven when $q$ is prime~\cite{BallMDS}. 

The following theorem by Higgs is known as the \emph{Scum Theorem}, and will be of importance later in this paper. It significantly restricts the sets $A\subseteq B\subseteq E$ that one must consider in order to find all minors of $M$ as $M|B/A$. %The following version of the scum theorem can be found in \cite{oxley11}.

\begin{theorem}[Proposition 3.3.7 in \cite{crapo70}] \label{thm:scum}
Let $N=(E_N,\rho_N)$ be a minor of a matroid $M(E_M,\rho_M)$. Then there is a pair of sets $A \subseteq B\subseteq E_M$ with $\rho_M(A) = \rho_M(E_M)-\rho_N(E_N)$ and $\rho_M(B)=\rho_M(E_M)$, such that $M|B/A \cong N$. Further, if $N$ has no loops, then $A$ can be chosen to be a flat of $M$. 
\end{theorem}

\subsection{Fundamentals on cyclic flats}

Before we define and give the properties of the cyclic flats, we need a minimal background on posets and lattices. We refer the reader to \cite{stanley2011} for further information about these objects. 

A \emph{partially ordered set} $P$ (or \emph{poset}, for short) is a set together with a partial order $\leq$. For $x,y \in P$, we say that $y$ \emph{covers} $x$ or $x$ is \emph{covered} by $y$, denoted by $x \lessdot y$, if $x \leq y$, $x \neq y$, and there is no $z \in P$ different from $x$ and $y$ such that $x \leq z \leq y$. An \emph{upper bound} of $x$ and $y$ is an element $u \in P$ satisfying $x \leq u$ and $y \leq u$. The \emph{join} of $x$ and $y$, denoted by $x \vee y$ if it exists, is the least upper bound. Dually, the \emph{meet} $x \wedge y$ is the greatest lower bound. If $P$ has an element $0_{P}$ such that $0_{P} \leq x$ for all $x \in P$, then $0_{P}$ is called the \emph{bottom element} of $P$. Similarly, if $P$ has an element $1_{P}$ such that $x \leq 1_{P}$ for all $x \in P$, then $1_{P}$ is called the \emph{top element} of $P$. 

A \emph{lattice} is a poset $L$ for which every pair of elements has a meet and join. It is not difficult to see that every finite lattice has a bottom element and top element. For a lattice $L$ and $x \in L$, then $x$ is an \emph{atom} of $L$ if $x$ covers $0_{L}$. A lattice $L$ is said to be \emph{atomic} if every element of $L$ is the join of atoms. Dually, $x \in L$ is a \emph{coatom} if $x \lessdot 1_{L}$ and $L$ is \emph{coatomic} if every element of $L$ is the meet of coatoms. 

The main tool from matroid theory in this paper are the cyclic flats. We will define them using the closure and cyclic operator.   

Let $M = (E, \rho)$ be a matroid. The \emph{closure} operator $\mathrm{cl}:2^E \rightarrow 2^E$ and \emph{cyclic} operator $\mathrm{cyc}: 2^E \rightarrow 2^E$ are defined by
\begin{align*}
\cl(X) &= \{e \in E : \rho(X \cup e) = \rho(X)\},\\
\cyc(X) &= \{e \in X : \rho(X - e) = \rho(X)\}.
\end{align*}
A subset $X \subseteq E$ is a \emph{flat} if $\cl(X) = X$ and a \emph{cyclic set} if $\cyc(X) = X$. Therefore, $X$ is a \emph{cyclic flat} if 
$$
\rho(X \cup y) > \rho(X) \quad \hbox{and} \quad \rho(X - x) = \rho(X)
$$
for all $y \in E-X$ and $x \in X$. The collection of flats, cyclic sets, and cyclic flats of $M$ are denoted by $\mathcal{F}(M)$, $\mathcal{U}(M)$, and $\mathcal{Z}(M)$, respectively. 

It is easy to verify, as in~\cite{bonin08}, that the closure operator induces flatness and preserves cyclicity, and that the cyclic operator induces cyclicity and preserves flatness. Thus we can write 
\begin{equation}\label{eq:clcyc}
\cl: \left\{\begin{split} 2^E & \to\FF(M)\\ \UU(M)&\to\ZZ(M),\end{split}\right.\quad\textrm{ and }\quad\cyc: \left\{\begin{split} 2^E & \to\UU(M)\\ \FF(M)&\to\ZZ(M).\end{split}\right.
\end{equation} 
In particular, for any set $X\subseteq E$, we have $\cyc(\cl(X)) \in\ZZ(M)$ and $\cl(\cyc(X)) \in\ZZ(M)$. Moreover, the closure and cyclic operators are order preserving in that $X\subseteq Y$ implies $\cyc(X)\subseteq\cyc(Y)$ and $\cl(X)\subseteq\cl(Y)$, so the maps in~\eqref{eq:clcyc} can be considered as order-preserving poset maps.
The following duality properties of flats, cyclic sets, and cyclic flats are easy to verify.

\begin{proposition}
Let $M = (E,\rho)$ be a matroid and $X, Y \subseteq E$.
\begin{enumerate}[(i)]
%\item $\cl(X) \in \FF(M)$.
%\item If $X \subseteq Y$, then  $\cl(X) \subseteq  \cl(Y)$.
%\item If $F \in \FF(M)$ and $X \subseteq F$, then $\cl(X) \subseteq F$.
%\item $\cyc(X) \in \UU(M)$.
%\item If $X \subseteq Y$, then $\cyc(X) \subseteq  \cyc(Y)$.
%\item If $U \in \UU(M)$ and $U \subseteq X$, then $U \subseteq \cyc(X)$.
%\item If $F \in \FF(M)$, then $\cyc(F) \in \ZZ(M)$.
%\item If $U \in \UU(M)$, then $\cl(U) \in \ZZ(M)$.
\item $\FF(M) = \{E - X : X \in \UU(M^*)\}$.
\item $\UU(M) = \{E - X : X \in \FF(M^*)\}$.
\item $\ZZ(M) = \{E - X : X \in \ZZ(M^*)\}$.
\end{enumerate}
\end{proposition}

Two basic properties of cyclic flats of a matroid are given in the following proposition.

\begin{proposition} [\cite{bonin08}] 
Let $M = (E,\rho)$ be a matroid and $\mathcal{Z}$ the collection of cyclic flats of $M$. Then, 
\begin{enumerate}[(i)]
\item $ \rho(X) = \min \{\rho(F) + | X \setminus F | : F \in \mathcal{Z}\} \hbox{, for } X \subseteq E,$
\item $(\mathcal{Z}, \subseteq)$ is a lattice with $X \vee Y = \cl(X \cup Y)$ and $X \wedge Y = \cyc( X \cap Y)$ for $X,Y \in \mathcal{Z}$
\item $1_{\mathcal{Z}}=\cyc(E)$ and $0_{\mathcal{Z}}=\cl(\emptyset)$.
\end{enumerate}
\end{proposition}

That $E$ together with the cyclic flats and their ranks together defines the matroid $M = (E,\rho)$ uniquely can be concluded from (i) in the proposition above. It is thus natural to cryptomorphically define matroids via an axiom scheme for their cyclic flats, as was done independently in~\cite{bonin08} and \cite{sims80}. This gives a compact way to represent and construct matroids. 

\begin{theorem} [\cite{bonin08} Th. 3.2 and \cite{sims80}] \label{th:Z-axiom}
Let $\mathcal{Z} \subseteq 2^E$ and let $\rho$ be a function $\rho: \mathcal{Z} \rightarrow \mathbb{Z}$. There is a matroid $M$ on $E$ for which $\mathcal{Z}$ is the set of cyclic flats and $\rho$ is the rank function restricted to the sets in $\mathcal{Z}$, if and only if
\begin{itemize}
\item[(Z0)] $\ZZ$ is a lattice under inclusion.
\item[(Z1)] $\rho(0_{\mathcal{Z}}) = 0$.
\item[(Z2)] If $X,Y \in \mathcal{Z}$ and $X \subsetneq Y$, then $0 < \rho(Y) - \rho(X) < | Y | - | X |$
\item[(Z3)] $\rho(X) + \rho(Y) \geq \rho(X \vee Y) + \rho(X \wedge Y) + | (X \cap Y) - (X \wedge Y)  |$ for all $X, Y\in \ZZ$.
\end{itemize}
\end{theorem}

\begin{definition}
A matroid is \emph{non-degenerate} if it does not have any loops or isthmuses. A matroid which has a loop or isthmus is \emph{degenerate}.
\end{definition}

For any matroid $M=(E,\rho)$, we observe that 
$$
0_\ZZ = \cl(\emptyset) = \{e \in E : \rho(e) = 0\} \hand 1_\ZZ = \cyc(E) = \{e \in E : \rho(E - e) = \rho(E)\}.
$$
Hence, $0_\ZZ$ and $E- 1_\ZZ$ are equal to the collection of loops and isthmuses, respectively. Consequently, we obtain the following lemma which gives a characterization of non-degenerate matroids via cyclic flats.

\begin{lemma} \label{lem:characterization_non-degenerate_via_Z}
A matroid is non-degenerate if and only if $0_\ZZ = \emptyset$ and $1_\ZZ = E$.
\end{lemma}

%The cyclic flats of a linear code $C$ of $\F^E$ can be described as sets $X \subseteq E$ such that 
%$
%\dim(C|(X \cup y)) > \dim(C|X)$ for all $y \in E -X$ and $\dim(C|(X - x)) = \dim(C|X)$
%for all $x \in X$. Thus, we have the following immediate proposition.

As an immediate consequence of the definition of cyclic flats, we have the following characterisation of uniform matroids by their cyclic flats.

\begin{proposition}\label{prop:uniform}
Let $M = (E, \rho)$ be a matroid and let $0<k<n$ be positive integers. The following are equivalent:
\begin{enumerate}[(i)]
\item $M$ is the uniform matroid $U_n^k$
\item $\ZZ = \ZZ(M)$ is the two element lattice with bottom element $0_{\ZZ}=\emptyset$, top element $1_{\ZZ} = E$ and $\rho(1_{\ZZ}) = k$
\end{enumerate}
\end{proposition}

Finally, it was proven in \cite{westerback15} that the cyclic flats determine the minimum distance of a non-degenerate matroid. 

\begin{proposition}[\cite{westerback15}]
Let $\rmatroid$ be a non-degenerate matroid. Then the minimum distance $d$ satisfies
\[
d=\eta(E)+1-\max\{ \eta(Z) : Z \in \ZZ(M)-E\}.
\]
\end{proposition}

 %Intro + Preliminaries
\section{Cyclic Flats of Minors}

In order to identify uniform minors of the matroid $M$, we will take the detour of identifying the cyclic flats of an arbitrary minor $M_{[A,B]}$. We will then use the fact that the minor in question is uniform if and only if $\ZZ(M_{[A,B]})=\{\emptyset, B-A\}$, as in Proposition~\ref{prop:uniform}. Our first interest is in  the case when $X$ and $Y$ are themselves cyclic flats. In this case we have a straightforward characterisation of cyclic flats in $M|Y/X$, via the following two lemmas:

\begin{lemma}\label{flats}
Let $M$ be a matroid, and let $X\subseteq Y\subseteq E(M)$ be two sets with $Y\in\FF(M)$. Then $\FF(M|Y/X)=\{F\subseteq Y-X , F\cup X\in\FF(M)\}$.
\end{lemma}

\begin{proof}
A set $S$ is flat in $M|Y/X$ precisely if $\rho(S\cup X\cup i)>\rho(S\cup X)$ for all $i\in (Y-X)- S$. Since $Y$ is a flat, the inequality $\rho(S\cup X\cup i)>\rho(S\cup X)$ will hold for all $i\in \bar{Y}$ regardless of $S$. Thus, $S$ is flat in $M|Y/X$ if and only if $S\cup X$ is flat in $M$.
\end{proof}

\begin{lemma}\label{cyclic}
Let $M$ be a matroid, and let $X\subseteq Y\subseteq E(M)$ be two sets with $X\in\UU(M)$. Then $\UU(M|Y/X)=\{U\subseteq Y-X , U\cup X\in\UU(M)\}$.
\end{lemma}
This is the dual statement, and thus an immediate consequence, of Lemma~\ref{flats}. We write out the proof explicitly only for illustration.
\begin{proof}
A set $S$ is cyclic in $M|Y/X$ precisely if $\rho((S\cup X)- i)=\rho(S\cup X)$ for all $i\in S$. For $i\in X$, this will hold regardless of $S$, since $X$ is cyclic. Thus, $S$ is cyclic in $M|Y/X$ if and only if $S\cup X$ is cyclic in $M$.
\end{proof}

The previous lemmas give the following immediate corollary:
\begin{corollary}\label{corr:interval}
Let $M=(E, \rho)$ be a matroid, and let $X\subseteq Y\subseteq E(M)$ be two sets with $X\in\UU(M)$ and $Y\in \FF(M)$. Then $\ZZ(M|Y/X)=\{Z\subseteq Y-X , Z\cup X\in\ZZ(M)\}$, with the rank function $\rho_{|Y/X}(Z)=\rho(Z\cup X)-\rho(X)$.
\end{corollary}

In particular, for any $X\subseteq Y\subseteq E$, $\ZZ(M|Y/X)$ is isomorphic to an interval in $\ZZ(M)$. As a consequence, we get a sufficient condition for uniformity of minors, that only depends on the Hasse diagram of $\ZZ(M)$. 

\begin{theorem}\label{thm:edge}
Let $X$ and $Y$ be two cyclic flats in $M$ with $X\lessdot_{\ZZ(M)} Y$. Let $n=|Y|-|X|$ and $k=\rho(Y)-\rho(X)$. Then $M|Y/X\cong U_n^k$.
\end{theorem}

\begin{proof}
By Corollary~\ref{corr:interval}, we have $\ZZ(M|Y/X)=\{\emptyset,Y-X\}$ as there are no cyclic flats with $X\subset Z\subset Y$. Again by Corollary~\ref{corr:interval}, we have $\rho(Y-X)=\rho(Y)-\rho(X)=k$. Thus, by Proposition~\ref{prop:uniform}, $M|Y/X$ is the uniform matroid $U_n^k$.
\end{proof}

\begin{corollary}\label{corr:hasse}
Let $M$ be a matroid that contains no $U_n^k$ minors. Then, for every edge $X\lessdot_{\ZZ(M)} Y$ in the Hasse diagram of $\ZZ(M)$, we have $\rho(Y)-\rho(X)<k$ or $\eta(Y)-\eta(X)<n-k$.
\end{corollary}
\begin{proof}
Assume for a contradiction that $X\lessdot_{\ZZ(M)} Y$ has $\rho(Y)-\rho(X)=k'\geq k$ and $\eta(Y)-\eta(X)=n'-k'\geq n-k$. Then by Theorem~\ref{thm:edge}, $M|Y/X\cong U_{n'}^{k'}$, and so contains $U_n^k$ as a minor by Lemma~\ref{lm:subunif}.
\end{proof}

Now, we are going to need formulas for how to compute the lattice operators in $\ZZ(M|Y/X)$ in terms of the corresponding operators in $\ZZ(M)$. These can be derived from corresponding formulas for the closure and cyclic operator. To derive these, we will need to generalize Corollary~\ref{corr:interval} to the setting where the restriction and contraction are not necessarily performed at cyclic flats.

\begin{theorem}\label{thm:cf_formulas}
For $X \subseteq Y \subseteq E$, we have
\begin{enumerate}
\item $\ZZ(M|Y) = \{ \cyc(Z \cap Y) : Z \in \ZZ(M) \}$.
\item $\ZZ(M/X) = \{ \cl(X \cup Z ) - X : Z \in \ZZ(M) \}$.
\item $\ZZ(M|Y/X) = \left\lbrace \cl\Big( X \cup \cyc \big( Z \cap Y\big) \Big) \cap \big(Y-X \big) : Z \in \ZZ(M) \right\rbrace \\
= \left\lbrace \cyc\Big( \cl \big(X \cup Z \big) \cap Y \Big) - X : Z \in \ZZ(M) \right\rbrace$.
\end{enumerate}
\end{theorem}

\begin{proof}
\begin{enumerate}
\item First, observe that the cyclic operator in $M|Y$ is the same as that in $M$, and that the flats in $M|Y$ are $\{ F \cap Y : F \in \FF(M) \}$. Thus we have \[\ZZ(M|Y)=\{\cyc(F\cap Y) : F\in\FF(M)\}\supseteq \{\cyc(Z\cap Y) : Z\in\ZZ(M)\}.\]  
On the other hand, let $A \in \ZZ(M|Y)$, so $\cyc(A) = A$ and $\cl(A) \cap Y = A$. But the closure operator preserves cyclicity, so $\cl(A)\in\ZZ(M)$. We then observe that \[A = \cyc(\cl(A) \cap Y)\in \{\cyc(Z\cap Y) : Z\in\ZZ(M)\}.\] This proves the reverse inclusion \[\ZZ(M|Y)=\{\cyc(F\cap Y) : F\in\FF(M)\}\subseteq \{\cyc(Z\cap Y) : Z\in\ZZ(M)\}.\]  

\item This is the dual statement of Statement 1, and so follows immediately by applying Statement 1 to the matroid $M^*|\bar{X}/\bar{Y}$. \item We first apply Statement 2 and then Statement 1 to the restricted matroid $M|Y$, and get
\[
\ZZ(M|Y/X) = \{ \cl_{|Y}(X \cup \cyc(Z \cap Y) ) - X : Z \in \ZZ(M) \}.
\]
Since $\cl_{|Y}(T) = \cl(T) \cap Y$ if $T \subseteq Y$, then we have
\[
\ZZ(M|Y/X) = \{ \cl(X \cup \cyc(Z \cap Y)) \cap (Y - X) : Z \in \ZZ(M) \}.
\]

For the second equality in Statement 3, we need to study the operator $\cyc_{/X}$. Suppose $T \subseteq E-X$. Using duality and the formula for $\cl_{| \bar{X}}$, we find that 
\[
\cyc_{/X}(T) = \cyc(X \cup T ) - X.
\]

%Details:
%\begin{align*}
%cyc_{M/X}(T) & = (E-X) - cl_{M^{\ast}\setminus X }( (E-X) -T ) \\
%& = (E-X) - (cl_{M^{\ast}} ((E-X)-T) \cap (E-X) ) \\
%& = ((E-X) - cl_{M^{\ast}} ( (E-X) - T ) ) \cup ((E-X) - (E-X) ) \\
%& = (E-X) - cl_{M^{\ast}} ( E- (X\cup T ) ) \\
%& = (E-X) - (E-cyc_{M}(X\cup T ) ) \\
%& = ((E-X) \cap cyc_{M}(X\cup T ) ) \cup ((E-X) - E ) \\
%& = cyc_{M}(X\cup T ) - X
%\end{align*}

Now we are ready to prove the last equality. Applying first Statement 1 and then Statement 2 to the contracted matroid $M/X$, we get
\[
\ZZ(M|Y/X) = \{ \cyc_{/X} ( ( \cl(X \cup Z ) - X ) \cap Y ) : Z \in \ZZ(M) \}.
\]
Applying the formula for $\cyc_{/X}$, we obtain
\begin{align*}
\ZZ(M|Y/X) & = \{ \cyc(( \cl(X \cup Z ) \cap Y-X ) \cup X ) - X : Z \in \ZZ(M) \}\\ & = \{ \cyc ( \cl( X\cup Z) \cap Y ) - X : Z \in \ZZ(M) \},
\end{align*}
where the last equality follows as $X \subseteq \cl(X \cup Z )$. This concludes the proof.
\end{enumerate}
\end{proof}

\section{Sufficient Conditions for Uniformity}

From Statement 3 of Theorem~\ref{thm:cf_formulas} we get two surjective maps $\ZZ(M)\to\ZZ(M|Y/X)$, given by $$Z\mapsto\cyc ( \cl( X\cup Z) \cap Y ) - X \hbox{ and }Z\mapsto \cl\Big( X \cup \cyc \big( Z \cap Y\big) \Big) \cap \big(Y-X \big)$$ respectively. To identify uniform minors in $M$, we need to detect $X$ and $Y$ such that either, and thus both, of these maps have image $\{\emptyset, Y-X\}$. %Using this, we will be able to find some conditions on the sets $X$ and $Y$, as well as on the matroid itself.

\subsection{Minors Given by Restriction or Contraction Only}

We will begin by considering a simpler case when the minor is the result of a restriction only, \emph{i.e.,} when the minor is given by $M|Y$. So, let $M=(E, \rho)$ be a matroid and $Y$ an arbitrary subset of $E$. We can use Corollary \ref{corr:interval} to restrict the amount of information we need to consider. Indeed, it is straightforward to see that
\begin{equation}
\label{minor_1}
 M|Y = M|\cl(Y) \setminus (\cl(Y)-Y).
\end{equation} 
Then, Theorem \ref{thm:cf_formulas} states that the cyclic flats of $M|Y$ depend only on the cyclic flats of $M|\cl(Y)$. Furthermore, according to Corollary \ref{corr:interval} the cyclic flats of $M|\cl(Y)$ are exactly the cyclic flats of $M$ contained in $\cl(Y)$. Hence, we can restrict the study to the case when $M=(E, \rho)$ is a matroid and $Y$ is a subset of full rank. Define $k:=\rho(Y)$ and $n:=|Y|$. With this setup, we obtain the following theorem.

\begin{theorem}
\label{thm_restriction_case}
Let $M=(E,\rho)$ be a matroid and $Y$ a subset of full rank. $M|Y$ is isomorphic to the uniform matroid $U_n^k$ if and only if either $Y$ is a basis of $M$ or the following two conditions are satisfied:
\begin{enumerate}
\item $Y$ is a cyclic set of $M$.
\item For all $Z \in \ZZ(M)$ with $\rho(Z)<k$, $Z \cap Y$ is independent in $M$.  
\end{enumerate}
\end{theorem}

Before stating the proof, we will need one useful lemma about the properties of the closure and cyclic operators. 

\begin{lemma}
\label{lemma:cl_cyc}
Let $M=(E,\rho)$ be a matroid and $Y \subseteq E$. Then
\begin{enumerate}
\item $\cl(\cyc(Y)) \cap Y = \cyc(Y)$.
\item $\cyc(\cl(Y)) \cup Y = \cl(Y)$.
\end{enumerate}
\end{lemma}

The proof of Lemma~\ref{lemma:cl_cyc} is straightfrorward from the definition of the operators together with the submodularity of the rank function. Details of the proof can be found in~\cite{bonin08}. 

\begin{corollary}
\label{cor:cl_cyc}
Let $M=(E,\rho)$ be a matroid and $Y \subseteq E$. Then
\begin{enumerate}
\item $Y\in\UU$ if and only if $\cl(\cyc(Y)) = \cl(Y)$.
\item $Y\in\FF$ if and only if $\cyc(\cl(Y)) = \cyc(Y)$.
\end{enumerate}
\end{corollary}
\begin{proof}
The right implications are immediate as $Y\in\UU$ means that $\cyc(Y) = Y$ and $Y\in\FF$ means that $\cl(Y) = Y$. For the left implication in the first statement, notice that if $\cl(\cyc(Y)) = \cl(Y)$, then by Lemma~\ref{lemma:cl_cyc} we have $$Y=\cl(Y)\cap Y=\cl(\cyc(Y))\cap Y=\cyc(Y)\cap Y=\cyc(Y),$$ so $Y$ is cyclic. The second statement is the dual of the first.
\end{proof}

We now present the proof of Theorem \ref{thm_restriction_case}. 

\begin{proof}
We know that $M|Y \cong U_n^k$ if and only if $\ZZ(M|Y) = \{ \emptyset, Y \}$. On the other hand, we know by Theorem \ref{thm:cf_formulas}, that 
\[
\ZZ(M|Y) = \{ \cyc(Z \cap Y ) : Z \in \ZZ(M) \}.
\]
Then we have that $M|Y \cong U_n^k$ if and only if  $\cyc(Z \cap Y) \in \{ \emptyset, Y \}$ for all $ Z \in \ZZ(M)$.

Now consider the cyclic flat $Z_{Y} := \cl(\cyc(Y))$. Using Lemma \ref{lemma:cl_cyc}, we have that 
\[
\cyc(Z_{Y} \cap Y) = \cyc(Y).
\]
Two cases can occur. If $\cyc(Y) = \emptyset$ then $Y$ was a basis of $M$ and we end up with a minor isomorphic to $U_k^k$. If not, then $\cyc(Y)$ must be equal to $Y$. So, we have that $Y \in \ZZ(M|Y)$ if and only if $Y$ is a cyclic set and we obtain the first condition. Since $Y$ already has full rank, there is only one cyclic flat that contains $Y$, namely $\cl(\cyc(Y)) = E$. Therefore, for every other cyclic flat $Z$, \emph{i.e.,} for all $Z$ with $\rho(Z)<k$, we have that
\[
\cyc(Z \cap Y) \subseteq Z \cap Y \neq Y.
\]

However, by Theorem \ref{thm:cf_formulas}, $\cyc(Z \cap Y ) $ is a cyclic flat of $M|Y$. Thus, for all $Z \in \ZZ(M)$ with $\rho(Z) <k$, we have $\cyc(Z\cap Y) = \emptyset$, or equivalently, $Z \cap Y$ is independent. Notice that, combined with the first condition, this implies immediately that $\ZZ(M|Y) = \{ \emptyset, Y \}$. This concludes the proof. 
\end{proof}

\begin{corollary}
Let $M=(E,\rho)$ be a matroid and $Y$ a subset of full rank. If $M|Y$ is isomorphic to $U_n^k$, then the ground set $E$ must be a cyclic flat, \emph{i.e.,} $E \in \ZZ(M)$. 
\end{corollary}

%\begin{example}
%
%We will look at the matroid $M_{C}$ arising from the binary matrix in Example \ref{ex:code_matroid}. Since it is a binary matroid, we cannot find a minor isomorphic to $U_4^2$. Using the lattice of cyclic flats of Fig. \ref{fig:LCF_Ex}, we can find a minor isomorphic to $U_3^2$ if we look at, \emph{e.g.},  $M_{C}|\{1,2,3\}$. We can also find a minor isomorphic to $U_4^3$. To this end, we look at $M_{C}|\{1,2,4,5\}$. The set $\{1,2,4,5\}$ is indeed a cyclic set because $\cyc(\{1,2,4,5 \} ) = \{1,2,4,5\}$. Furthermore, there is no cyclic flat properly contained in $\{1,2,4,5\}$, and every intersection of a cyclic flat different from $E$ with the set $\{1,2,4,5 \}$ gives us an independent set. Thus, the conditions from the previous theorem are met and we get a minor isomorphic to $U_4^3$. 
%
%\end{example}

Now we can do the same for $M/X$ and use duality to get back to the restriction case.  By minor properties, we have
\begin{equation}
\label{minor_2}
 M/X = M/\cyc(X) / (X-\cyc(X)).
\end{equation}
 
Then, Corollary \ref{corr:interval} states that the cyclic flats of $M/\cyc(X)$ are the cyclic flats of $M$ that contain $\cyc(X)$. Thus, we will consider a matroid $M=(E,\rho)$ and $X$ an independent subset of $E$. Define $k:=\rho(E)-\rho(X)$ and $n:=|E-X|$. Then, we have the following dual statement of Theorem \ref{thm_restriction_case}. 

\begin{theorem}
Let $M=(E,\rho)$ be a matroid and $X$ an independent subset of $E$. $M/X$ is isomorphic to the uniform matroid $U_n^k$ if and only if either $X$ is a basis of $M$ or the following two conditions are satisfied:
\begin{enumerate}
\item $X$ is a flat of $M$.
\item For all $Z \in \ZZ(M)-0_{\ZZ}$, we have $\cl(X \cup Z ) = E$.
\end{enumerate}
\end{theorem}

\begin{corollary}
Let $M=(E,\rho)$ be a matroid and $X$ an independent subset of $E$. If $M/X$ is isomorphic to $U_n^k$, then the empty set must be a cyclic flat, \emph{i.e.,} $\emptyset \in \ZZ(M)$.
\end{corollary}

\subsection{Minors Given by Both Restriction and Contraction}

This part combines the two previous situations into a more general statement. We will see that, when we allow both a restriction and a contraction to occur, we lose some conditions on the matroid that are then replaced by conditions on the sets used in the minor. 
%\textcolor{red}{we mentioned before (CAMINOTE: where, can we give a precise ref?)}

Let $M=(E,\rho)$ be a matroid and $X \subset Y \subseteq E$ two sets. Combining minor properties (\ref{minor_1}) and (\ref{minor_2}) and Corollary \ref{corr:interval}, it is sufficient to only consider the cyclic flats  between $\cyc(X)$ and $\cl(Y)$. In addition, we want to avoid some known cases, namely when $Y$ is independent (we will obtain $U_n^n$) and when $X$ has full rank (we will obtain $U_n^0$). Define $k:=\rho(E)-\rho(X)$ and $n:=|Y-X|$. We get the following theorem.  

\begin{theorem}
Let $M=(E,\rho)$ be a matroid and $X \subset Y \subseteq E$ two sets such that $Y$ is a dependent full-rank set and $X$ is an independent set with $\rho(X)<\rho(E)$. The minor $M|Y/X$ is isomorphic to a uniform matroid $U_n^k$ if and only if
\begin{enumerate}
\item $\cl(X) \cap Y = X$,
\item $Y-X \subseteq \cyc(Y)$, and
\item for all $Z \in \ZZ(M)$ either $Z \cap Y$ is independent or\\ $\cl(X \cup \cyc(Z \cap Y)) = E$.
\end{enumerate}
\end{theorem}

\begin{proof}
Using Theorem \ref{thm:cf_formulas}, we have that $M|Y/X \cong U_n^k$ if and only if 
\[
\cl(X \cup \cyc(Z \cap Y) ) \cap (Y-X) \in \{ \emptyset, Y-X \} \text{ for all } Z \in \ZZ(M).
\]
In particular, it holds for $Z = 0_{\ZZ}$. Let $Z'_{0} := \cl(X \cup \cyc(0_{\ZZ} \cap Y))$. Using the properties of the closure, we have
\[
\cl(X) \subseteq Z'_{0} \subseteq \cl(X \cup 0_{\ZZ} ) = \cl(X).
\]
Thus, we have a chain of equalities and, in particular, $\rho(Z'_{0}) < \rho(E)$. This means that $Z'_{0} = \emptyset$ in order to have $\emptyset \in \ZZ(M|Y/X)$. But, since $Z'_{0} = \cl(X)$ then $Z'_{0} = \emptyset$ is equivalent to  $\cl(X) \cap Y = X$ and Condition~1 is proved.

Now, consider the cyclic flat $Z_{Y} := \cl(\cyc(Y))$. First, using again Lemma \ref{lemma:cl_cyc}, we have $\cyc(Z_{Y} \cap Y ) = \cyc(Y)$. Since $Y$ is a dependent subset, $\cyc(Y) \neq \emptyset$ and thus $X \subsetneq X \cup \cyc(Y)$. Then, the closure cannot be contained in $X$ and we must have
\[
\cl(X \cup \cyc(Y) ) \cap (Y-X) = Y-X.
\]

Define $Z'_{Y} := \cl(X \cup \cyc(Y) )$. The above equality implies that $Y-X \subseteq Z'_{Y}$. On the other hand, we have that $X \subseteq Z'_{Y}$. Then, $Y \subseteq Z'_{Y}$ and $Z'_{Y} = E$. In particular, we must have $Y-X \subseteq \cyc(Y)$. Indeed, assume by contradiction that there exists $a \in Y-X$ and $a \notin \cyc(Y)$. Then, by definition of the cyclic operator, $\rho(Y-a)<\rho(Y)$. But since $a \notin \cyc(Y)$ and $a \notin X$ then $X \cup \cyc(Y) \subseteq Y-a$. This implies that 
\[
\rho(X \cup \cyc(Y)) \leq \rho(Y-a) < \rho(Y),
\]
which is a contradiction. The condition $Y-X \subseteq \cyc(Y)$ is also sufficient to guarantee that $Y-X \in \ZZ(M|Y/X)$. This proves Condition 2.  

Finally, for every other cyclic flat of $M$, $\cl(X \cup \cyc(Z \cap Y) ) \cap (Y-X) = \emptyset$ if and only if $\cyc(Z \cap Y) \subseteq X$. But since $X$ is independent, this is equivalent to $Z \cap Y$ being independent. On the other hand, 
\[
\cl(X \cup \cyc(Z \cap Y) ) \cap (Y-X) = Y-X \text{ if and only if } \cl(X \cup \cyc(Z \cap Y) ) = E.
\]
This concludes the proof. 
\end{proof}

This theorem and the proof are only based on the first representation of the cyclic flats of $M|Y/X$ in Theorem \ref{thm:cf_formulas}. We can also state an equivalent theorem obtained using the second representation in Theorem \ref{thm:cf_formulas}.

\begin{theorem}
Let $M=(E, \rho)$ be a matroid and $X \subset Y \subseteq E$ two sets such that $Y$ is a dependent full-rank set and $X$ is an independent set with $\rho(X) < \rho(E)$. The minor $M|Y/X$ is isomorphic to a uniform matroid $U_{k,n}$ if and only if
\begin{enumerate}
\item $\cl(X) \cap Y = X$,
\item $Y-X \subseteq \cyc(Y)$, and
\item for all $Z \in \ZZ(M)$ either $\cl(X \cup Z) \cap Y$ is independent or \\ $\cl(X \cup Z ) = E$.
\end{enumerate}
\end{theorem}

We conclude this section with a lemma that will be used later in our analysis, about pairs $Z\subset Z_1$ of cyclic flats with rank difference equal to one. 

\begin{lemma}~\label{lm:rho1}
Let $Z, Z_1, Z_2$ be cyclic flats with $Z\subseteq Z_1$, $Z\subseteq Z_2$, $Z_1\not\subseteq Z_2$ and $\rho(Z_1)=\rho(Z)+1$. Then $Z_1\cap Z_2=Z$.
\end{lemma}

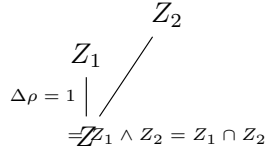
\begin{figure}
\centering
\resizebox{0.3\textwidth}{!}{%
\begin{tikzpicture}
\node[] (Z) at (0,0) {$Z$};  
\node[] (Z1) at (0,1) {$Z_1$};
\node[] (Z2) at (1,1.5) {$Z_2$};
\path [-] (Z) edge node[left] {\tiny{$\Delta\rho=1$}} (Z1);
\path [-] (Z) edge (Z2);
\node[] (t) at (1,0) {\tiny{$=Z_1\wedge Z_2=Z_1\cap Z_2$}};
\end{tikzpicture}
}
\caption{Illustration of Lemma~\ref{lm:rho1}.}
\label{fig:Intervals_Lemma6}
\end{figure}
\begin{proof}
The intersection $Z_1\cap Z_2$ of two flats is a flat of rank $<\rho(Z_1)$, which contains $Z$ by assumption. But $Z$ is a flat of rank $\rho(Z_1)-1$, so any set properly containing it has rank $\geq \rho(Z_1)$. It follows that $Z_1\cap Z_2=Z$.
\end{proof}

 % Z(M|B/A) from Z(M) + Uniform minors
%Input "u2n"
%Contains : F(M) from Z(M) + Char of U_2^n-avoiding from Z(M)

\section{Reconstructing the Lattice of Flats}
As the lattice of cyclic flats together with the induced rank function uniquely determines a matroid, it clearly also defines the lattice of flats $\FF(M)$. However, reconstructing $\FF(M)$ from $\ZZ(M)$ is not entirely straightforward. In order to do this, we will use the following notation.

\begin{definition}
Let $M=(E,\rho)$ be a matroid and let $A\subseteq E$. We will denote by $A^\FF$ the set $\{e\in E-A | A\cup\{e\}\in\FF(M)\}$. 
\end{definition}

In particular, we see that if $A$ is not a flat, then $A^\FF$ is at most a singleton, because otherwise $A$ could be written as an intersection $$A=\cap_{e\in A^\FF}(A\cup\{e\})$$ of flats, and would thus be a flat itself. On the other hand, if $A$ is a flat, then we get the following equivalent description of $A^\FF$.

\begin{lemma}\label{lm:AFrank}
Let $A$ be a flat in $M=(E,\rho)$ and $e\in E-A$. Then $e\in A^\FF$ if and only if  $\rho(A\cup\{e,f\})=\rho(A)+2$ for all $f\in E-A-\{e\}$.
\end{lemma}
\begin{proof}
As $A$ is a flat, $\rho(A\cup\{e\})=\rho(A)+1$ for all $e\in E-A$. For such $e$, $e\in A^\FF$ if and only if $A\cup \{e\}$ is a flat, which it is precisely if \[\rho(A\cup \{e\}\cup \{f\})>\rho(A\cup\{e\})=\rho(A)+1\] for all $f\in E-A-\{e\}$.
\end{proof}

From Lemma~\ref{lm:AFrank}, it easily follows that $A\mapsto A^\FF$ is an order-reversing set-valued map on the lattice of flats, and as a consequence also on the lattice of cyclic flats.

\begin{lemma}\label{lm:AFreverse}
Let $A$ and $B$ be flats in $M=(E,\rho)$ with $A\subseteq B$. Then $A^\FF\supseteq B^\FF$.
\end{lemma}
\begin{proof}
If $e\in B^\FF$, then by Lemma~\ref{lm:AFrank} we have $\rho(B\cup\{e,f\})=\rho(B)+2$ for all $f\in E-B-\{e\}$. By submodularity of $\rho$, we then also get $\rho(A\cup\{e,f\})=\rho(A)+2$ for all $f\in E-B-\{e\}$. Moreover, for $f\in B-A$, we have 
\begin{align*}
\rho(A\cup\{e,f\})-\rho(A)&= \left(\rho(A\cup\{e,f\})-\rho(A\cup \{f\})\right)+\left(\rho(A\cup \{f\})-\rho(A)\right) \\
&\geq \left(\rho(B\cup\{e\})-\rho(B)\right)+\left(\rho(A\cup\{f\})-\rho(A)\right)\\
&=1+1=2.
\end{align*}
Thus we have $\rho(A\cup\{e,f\})=\rho(A)+2$ for all $f\in E-A-\{e\}$, so $e\in A^\FF$.
\end{proof}

To an arbitrary $A\subseteq E$, we will now associate two intervals in $\ZZ(M)$, both of which will yield the singleton $\{A\}$ if $A$ is already a cyclic flat.

\begin{definition}\label{def:updown}
Let $A\subseteq E$, and let \[A^{\vee}=\bigvee_{\substack{Z\subseteq A\\ Z\in\ZZ(M)}}Z \,\mbox{ and } \, A^{\wedge}=\bigwedge_{\substack{Z\supseteq A\\ Z\in\ZZ(M)}}Z.\] 
\end{definition}

\begin{lemma}\label{lm:ZAorder}
Let $A\subseteq E$. Then we have the inclusions \[A^\vee\subseteq \cl(\cyc(A))\subseteq \cyc(\cl(A))\subseteq A^\wedge.\]
\end{lemma}
\begin{proof}
Every cyclic flat that is a subset of $A$ is also a subset of $\cyc(A)$, so $\bigcup_{\substack{Z\subseteq A\\ Z\in\ZZ(M)}}Z\subseteq \cyc(A)$. It then follows that \[A^{\wedge}=\cl\left(\bigcup_{\substack{Z\subseteq A\\ Z\in\ZZ(M)}}Z\right)\subseteq\cl(\cyc(A)).\] Dually, every cyclic flat that contains $A$ also contains $\cl(A)$, so $\bigcap_{\substack{Z\supseteq A\\ Z\in\ZZ(M)}}Z\supseteq \cl(A)$. It then follows that \[A^{\vee}=\cyc\left(\bigcap_{\substack{Z\supseteq A\\ Z\in\ZZ(M)}}Z\right)\supseteq\cyc(\cl(A)).\] Finally, as $A\subseteq \cl(A)$ we have that $\cyc(A)\subseteq \cyc(\cl(A))$, and as the latter is a flat by \eqref{eq:clcyc}, we have $\cl(\cyc(A))\subseteq \cyc(\cl(A))$.

\end{proof}

By Lemma~\ref{lm:ZAorder} we can define the intervals \[\ZZ(A)= [\cl(\cyc(A)), \cyc(\cl(A))]_{\ZZ(M)} \, \mbox{ and } \,\ZZ'(A)=[A^\vee , A^\wedge]_{\ZZ(M)},\] and observe that $\ZZ(A)\subseteq \ZZ'(A)$. 

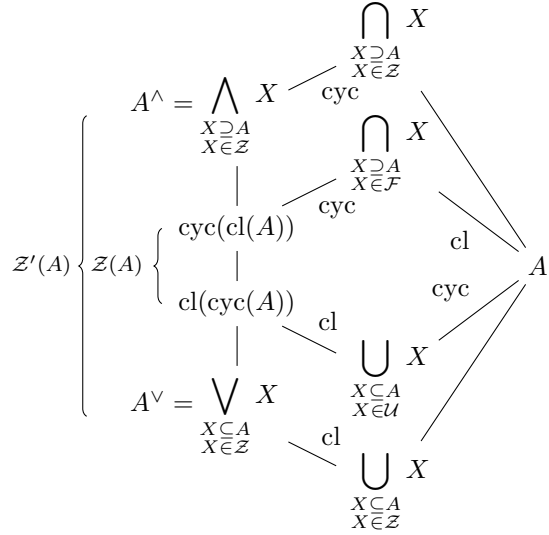
\begin{figure}
\centering
\resizebox{0.6\textwidth}{!}{
\begin{tikzpicture}
\node[] (A) at (2,0) {$A$};  
\node[] (top) at (0,3) {$\begin{aligned}\bigcap_{\substack{X\supseteq A \\ X\in \ZZ}}X\end{aligned}$};
\node[] (clA) at (0,1.5) {$\begin{aligned}\bigcap_{\substack{X\supseteq A \\ X\in \FF}}X\end{aligned}$};
\node[] (cycA) at (0,-1.5) {$\begin{aligned}\bigcup_{\substack{X\subseteq A \\ X\in \UU}}X\end{aligned}$};  
\node[] (bottom) at (0,-3) {$\begin{aligned}\bigcup_{\substack{X\subseteq A \\ X\in \ZZ}}X\end{aligned}$};  
\path [-] (A) edge (top);
\path [-] (A) edge node[below left] {$\cl$} (clA);
\path [-] (A) edge node[above left] {$\cyc$} (cycA);
\path [-] (A) edge (bottom);
\node[] (Ztop) at (-2,2) {$\begin{aligned}\bigwedge_{\substack{X\supseteq A \\ X\in \ZZ}}X\end{aligned}$};
\node[] (Awedge) at (-3,2.2) {$A^\wedge=$};
\node[] (ZclA) at (-2,0.5) {$\cyc(\cl(A))$};  
\node[] (ZcycA) at (-2,-0.5) {$\cl(\cyc(A))$};
\node[] (Zbottom) at (-2,-2) {$\begin{aligned}\bigvee_{\substack{X\subseteq A \\ X\in \ZZ}}X\end{aligned}$};  
\node[] (Avee) at (-3,-1.8) {$A^\vee=$};
\path [-] (Ztop) edge node[below right] {$\cyc$}  (top);
\path [-] (ZclA) edge node[below right] {$\cyc$} (clA);
\path [-] (ZcycA) edge node[above right] {$\cl$} (cycA);
\path [-] (Zbottom) edge node[above right] {$\cl$}  (bottom);
\path[-] (Ztop) edge (ZclA) ;
\path[-] (ZclA) edge (ZcycA);
\path[-] (ZcycA) edge (Zbottom);
\draw [decorate,decoration={brace,amplitude=3pt}]
(-3,-0.5) -- (-3,0.5) node [black,midway,xshift=-0.6cm] 
{\footnotesize $\ZZ(A)$};
\draw [decorate,decoration={brace,amplitude=3pt}]
(-4,-2) -- (-4,2) node [black,midway,xshift=-0.6cm] 
{\footnotesize $\ZZ'(A)$};
\end{tikzpicture}
}
\caption{Illustration of Definition~\ref{def:updown} and Lemma~\ref{lm:ZAorder}}
\label{fig:Intervals}
\end{figure}

The following lemma first occurred in~\cite{shoda16}. We state it here for completeness.
\begin{lemma}\label{lm:eqlattice}
Let $Z\in \ZZ$ satisfy $\rho(Z)+|A-Z|=\rho(A)$. Then $Z\in\ZZ(A)$.
\end{lemma}
\begin{proof}
Note that, for any $Z\subseteq E$, \[\rho(A)\leq \rho(A\cap Z) + |A-Z|\leq \rho(Z)+|A-Z|.\] The first inequality is satisfied with equality if and only if $\cyc(A)\subseteq A\cap Z$, which implies $\cyc(A)\subseteq Z$. We claim that the second inequality is satisfied with equality only if $Z\subseteq \cl(A)$. Indeed, if $Z\not\subseteq \cl(A)$ choose $z\in Z-\cl(A)$. Then $\rho(A)<\rho(A\cup\{z\})$, so by submodularity we have $\rho(Z\cap A)<\rho((Z\cap A)\cup\{z\})\leq\rho(Z)$. 

Therefore, any $Z$ satisfying $\rho(Z)+|A-Z|=\rho(A)$ must also satisfy $\cyc(A)\subseteq Z\subseteq\cl(A)$, so $\cl(\cyc(A))\subseteq\cl(Z)$ and $\cyc(Z)\subseteq\cyc(\cl(A))$. But, if $Z$ is a cyclic flat, then $Z=\cl(Z)=\cyc(Z)$, and it follows that \[Z\in[\cl(\cyc(A)),\cyc(\cl(A))]=\ZZ(A).\]
\end{proof}

\begin{lemma}\label{lm:flatcollapse}
Let $M=(E,\rho)$ be a matroid, and let $F\in\FF(M)$. Then $F^\vee=\cl(\cyc(F))=\cyc(\cl(F))$.
\end{lemma}
\begin{proof}
Since $F$ is a flat, we have $\cyc(\cl(F))=\cyc(F)\subseteq F$. Since $F^\vee$ contains all cyclic flats that are contained in $F$, it thus also contains $\cyc(\cl(F))$. The inclusions \[\cyc(\cl(F))\subseteq F^\vee\subseteq\cl(\cyc(F))\subseteq\cyc(\cl(F))\] now show that the sets are equal.
\end{proof}

\begin{corollary}\label{cor:flatcollapse}
Let $M=(E,\rho)$ be a matroid, and let $F\in\FF(M)$. Then $F^\vee=\cyc(F)\subseteq F$.\end{corollary}

We are now ready to present a result which explicitly reconstructs $\FF(M)$ from $\ZZ(M)$. Clearly, this result can immediately be dualized to obtain a description of $\UU(M)$.

\begin{proposition}~\label{prop:flatconditions}
Let $M=(E,\rho)$ be a matroid, and let $F\subseteq E$ be a set with $F^\vee\subseteq F$. Then the following are equivalent\begin{enumerate}[(i)]
\item\label{cond1} $F$ is a flat.
\item\label{cond3}$F^\vee= 1_{\ZZ(F)}$.
\item\label{cond2} For every $Z\in\ZZ'(F)-\{F^\vee\}$ it holds that \[|F\cap Z|-\rho(Z)< \eta(F^\vee).\]
\item\label{cond4} Every set $B$ with $F^\vee\subseteq B\subseteq F$ is a flat, and if $F\subsetneq F^\wedge$ then $|F|-\rho(F^\wedge)<\eta(F^\vee)$.
\item\label{cond5} For every set $B$ with $F^\vee\subseteq B\subseteq F$ and $B\subsetneq B^\wedge$ it holds that $|B|-\rho(B^\wedge)<\eta(B^\vee)$.
\end{enumerate}
\end{proposition}
\begin{proof} %First notice that if $F^\vee$ is a cyclic set contained in $F$, we also have $F^\vee\subseteq \cyc(F)$. We are now ready to prove the chain of inclusions.

\underline{\eqref{cond1}$\Rightarrow$\eqref{cond3}:} Assume $F\in\FF$, so $\cyc(F)\in \ZZ$. Then $\cyc(F)$ is a largest element of $\{Z\in \ZZ : Z\subseteq F\}$, so $\cyc(F)=F^\vee$. Moreover, we have $\cl(F)=F$, so $$1_{\ZZ(F)}=\cyc(\cl(F))=\cyc(F)=F^\vee.$$

\underline{\eqref{cond3}$\Rightarrow$\eqref{cond1}:}
As $F^\vee$ is a cyclic set contained in $F$, we also have $F^\vee\subseteq \cyc(F)$, so if \eqref{cond3} holds, then we have $$\cyc(\cl(F))=F^\vee\subseteq\cyc(F)\subseteq\cyc(\cl(F)).$$ Thus $\cyc(F)=\cyc(\cl(F)$, so $$\eta(F)=\eta(\cyc F)=\eta(\cyc(\cl(F))=\eta(\cl(F)).$$ As we also have $\rho(F)=\rho(\cl(F))$, it follows that $F=\cl(F)$, so $F$ is a flat.

\underline{\eqref{cond3}$\Rightarrow$\eqref{cond2}:}
Assume that $F^\vee= 1_{\ZZ(F)}$, so $\ZZ(F)=\{F^\vee\}$. For every $Z\in\ZZ'(F)-\{F^\vee\}$, since $F^\vee\subseteq F$, we then have $$|F|-\eta(F^\vee)=\rho(F^\vee)+ |F-F^\vee|=\rho(F)<\rho(Z)+|F-Z|=\rho(Z)+|F|-|F\cap Z|$$ by Lemma~\ref{lm:eqlattice}. Rewriting this equation, we immediately get 
$$|F\cap Z|-\rho(Z)< \eta(F^\vee).$$

\underline{\eqref{cond2}$\Rightarrow$\eqref{cond3}:}
Assume that \eqref{cond2} holds, and let $Z\in\ZZ'(A)-\{A^\vee\}$. Then we have \[\eta(F^\vee)-|F^\vee-F| = \eta(F^\vee) > |F\cap Z|-\rho(Z) = \eta(Z)-|Z-F|.\]
In particular, we see that $Z\neq 1_{\ZZ}$, so we must have
$A^\vee= 1_{\ZZ(A)}$.

\underline{\eqref{cond1} $\Rightarrow$\eqref{cond4} (assuming \eqref{cond3}$\Rightarrow$\eqref{cond1}$\Rightarrow$\eqref{cond2}):} 
Let $F$ be a flat with $F^\vee\subseteq B\subseteq F$. Then $F^\vee$ is a largest element of $\{Z\in \ZZ : Z\subseteq B\}$, so $B^\vee=F^\vee$. Also, $$B^\vee\subseteq 1_{\ZZ(BB)}\subseteq 1_{\ZZ(FF)}=F^\vee=B^\vee,$$ so we have equality $B^\vee= 1_{\ZZ(BB)}$. By the implication  \eqref{cond3}$\Rightarrow$\eqref{cond1}, $B$ is a flat. Now assume $F\subsetneq F^\wedge$. Since $F^\wedge\in\ZZ'(A)-\{A^\vee\}$, the implication \eqref{cond1}$\Rightarrow$\eqref{cond2} yields $$|F|-\rho(F^\wedge)=|F\cap F^\wedge|-\rho(F^\wedge)<\eta(F^\vee).$$

\underline{\eqref{cond5}$\Rightarrow$\eqref{cond2}:}
Let $Z\in\ZZ'(F)-\{F^\vee\}$ and set $B=F\cap Z$. 
Then $F^\vee\subseteq B\subseteq F$, and $1_{\ZZ(B)}\subseteq B^\wedge\subseteq Z$. If $1_{\ZZ(B)}\subsetneq Z$, then we get $$|B|-\rho(Z)=|B|-(\rho(Z)+|B-Z|)< |B|-\rho(B)=\eta(B)=\eta(\cyc(B))=\eta(F^\vee).$$
Now assume $1_{\ZZ(B)}= B^\wedge= Z$, so $B\subseteq Z= B^\wedge$.
If $B\subsetneq B^\wedge$, then by \eqref{cond5} we have $$|F\cap Z|-\rho(Z) = |B|-\rho(Z)\leq |B|-\rho(B^\wedge)<\eta(B^\vee)=\eta(F^\vee).$$ 
Finally, if $B=B^\wedge$, then $B$ is a cyclic flat, so  $B=F^\vee$ and $\ZZ(B)=\{B\}$. As $Z\neq F^\vee$ it follows that $Z\not\in\ZZ(B)$, so $$|B|-\rho(Z)=|B|-(\rho(Z)+|B-Z|)< |B|-\rho(B)=\eta(B)=\eta(\cyc(B))=\eta(F^\vee).$$

\underline{\eqref{cond4}$\Leftrightarrow$\eqref{cond5} (assuming \eqref{cond1}$\Leftrightarrow$\eqref{cond4}):}
This follows immediately by induction over the size of the set $B-F^\vee$.
\end{proof}

\section{Characterization of $U_{n}^{2}$ Avoiding Matroids from $\ZZ(M)$}
We will use the derived description of $\FF(M)$ together with Theorem~\ref{thm:scum} to detect whether $U_n^2$ is a minor of a matroid $M$ described by its cyclic flats. We use the notation $(F,E)_{\FF(M)}$ to denote the open interval between $F$ and $E$ in the lattice of flats $\FF(M)$. 

\begin{lemma} \label{lem:uniform_rank_2_minors_via_number_of_flats}
Let $M=(E,\rho)$ be a matroid, and $F\in\FF(M)$ be a flat with $\rho(E) - \rho(F) = 2$. Then $U_{n}^2$ is a minor of $M/F$ if and only if $|(F,E)_{\FF(M)}| \geq n$.
\end{lemma}

\begin{proof}
As $F$ is a flat, every proper superset of $F$ has rank $>\rho(E)-2$, so the minor $M|B/A$ has rank $<2$ whenever $F\subsetneq A$. Thus it is enough to show that there is a set $F\subseteq B\subseteq E$ with $M|B/F\cong U_n^2$ if and only if  $|(F,E)_{\FF(M)}| \geq n$. 

For the right implication, assume that there is a subset $F \subseteq B \subseteq E$ such that $\cl(A) = E$ and $M|B/F\cong U_{n}^2$. Further, let $(B-F) = \{b_1, \ldots, b_n \}$. %Then the nontrivial flats of $M|A/F_1$ are the singletons $\{b_1\}, \ldots, \{b_n\}$. 
Now $\cl(F \cup b_i)\in (F,E)_{\FF(M)}$ for $i = 1, \ldots,n$. Moreover, if $i\neq j$, as $\rho_{|B/F}(\{b_i,b_j\})=2$, we get that $$\rho(F\cup\{b_i,b_j\})=\rho(F)+2>\rho (F \cup b_i),$$ so $b(j)\not\in \cl(F \cup b_i)$. It follows that the flats $\cl(F \cup b_i)$ are all distinct, so $|(F,E)_{\FF(M)}| \geq n$.

For the left implication, let $B_1, \ldots, B_n$ be any $n$ different flats in $(F,E)_{\FF(M)}$. As $\rho(B_i) = \rho(F_1) + 1$, we get $B_i \cap B_j = F_1$ if $i \neq j$. Now, let $b_i \in B_i - F$ and let $B = F_1 \cup \{b_1, \ldots, b_n\}$. Then we get $\rho(F\cup\{b_i\})=\rho(F)+1$ for $i = 1, \ldots,n$, and $\rho(F\cup\{b_i,b_j\})=\rho(F)+2$ whenever $i\neq j$, and  $\rho(F\cup B)\leq\rho(E)=\rho(F)+2$. Thus we have $\rho_{|B/F}(A)=\min\{|A|, 2\}$,  
so $M|B/F \cong U_{2,n}$. 
\end{proof}

\begin{lemma}~\label{lm:u2n}
Let $M=(E,\rho)$ be a matroid. Then for $n \geq 3$, $U_n^2$ is a minor of $M$ if and only if there is a flat $F \in \FF(M)$ with $F\subseteq 1_{\ZZ}$ such that $\rho(1_\ZZ) - \rho(F) = 2$ and $|(F, 1_\ZZ)_\FF | \geq n$.
\end{lemma}
\begin{proof} If $e\in E$ is an isthmus, it is also an isthmus in $M|B/A$ for every $A\subseteq B\subseteq E$ with $e\in B$. Thus, as $U_n^2$ is non-degenerate, it can only occur as a minor $M|B/A$ where $B$ contains no isthmuses, so $B\subseteq 1_\ZZ$. By Theorem~\ref{thm:scum}, if $U_n^2$ is a minor of $M|{1_\ZZ}$, it is isomorphic to $M|B/F$ for some $F\subseteq B\subseteq 1_\ZZ$ with $\rho(F)+2=\rho(B)=\rho(E)$. But by Lemma~\ref{lem:uniform_rank_2_minors_via_number_of_flats}, this is equivalent to the condition $|(F, 1_\ZZ)_\FF | \geq n$.
\end{proof}

Note that the flats of corank $2$ in $1_\ZZ$ are easily identifiable via Proposition~\ref{prop:flatconditions}. If such a flat $F$ satisfies $F\cup e\in\FF$ for all $e\in 1_\ZZ-F$, then $M|1_\ZZ/F$ is uniform of rank $n$. We will therefore focus on flats $F$ of corank $2$ such that $1_\ZZ-F^\FF\neq \emptyset$. The key to detecting copies of $U_n^2$ from $\ZZ(M)$ now lies in determining $|(F, 1_\ZZ)_\FF |$ for such flats $F$, via studying certain antichains in $\ZZ(M)$. These antichains are defined next.

\begin{definition}\label{def:upsilon}
Let $M=(E, \rho)$ be a matroid and let $F\in\FF(M)$ be a flat with $\rho(F)=\rho(1_\ZZ)-2$. We define $\Upsilon(F)=\Upsilon_M(F)$ as the collection $$\Upsilon_M(F)=\left\{X\in\ZZ(M) : \rho(1_{\ZZ(M)})-1 = \rho(X) + |F-X| \right\},$$ and  $\bar{\Upsilon}(F)=\bar{\Upsilon}_M(F)$ as the collection of inclusion-maximal elements in $\Upsilon(F)$.
\end{definition}
%\ragnar{I skip one step in Thomas's definition, as it seems that the entire set of cyclic flats satisfying the equation is not relevant, but only the inclusion maximal ones.}

\begin{lemma}\label{lm:FH}
Let $M=(E, \rho)$ be a matroid, and let $F$ and $H$ be flats in $M$ with $F\lessdot_\FF H$ and $|H-F|\geq 2$. Then $H-F\subseteq\cyc(H)$.
\end{lemma}
\begin{proof}
Let $f\in H-F$. Then $F\subsetneq H-f$, so since $F$ is a flat we get $$\rho(F)<\rho(H-f)\leq \rho(H)=\rho(F)+1.$$ It follows that $\rho(H-f)= \rho(H)$ so $f\in\cyc (H)$.
\end{proof}

\begin{proposition}\label{lm:upsilon}
Let $M=(E, \rho)$ be a matroid and let $F\subseteq 1_{\ZZ}$ be a flat with $\rho(F)=\rho(1_\ZZ)-2$ and $1_\ZZ- F^\FF\neq\emptyset$. Then $$\bar{\Upsilon}(F)=\{\cyc(\cl(F\cup e)) : e\in 1_\ZZ- F^\FF\}.$$
\end{proposition}
\begin{proof} For $e\in 1_\ZZ - F^\FF$, let $H_e=\cl(F\cup e)$. Note that $F\lessdot_\FF H_e$ and $|H_e-F|\geq 2$. The proof will proceed in three steps: the first of two serve to show that $\cyc H_e\in \bar{\Upsilon}(F)$, and the third step shows that any $K\in\bar{\Upsilon}(F)$ can be written as $\cyc H_e$ for some $e\in 1_\ZZ- F^\FF$.

\underline{$\cyc(H_e)\in \Upsilon(F)$:}
We have 
$$\rho(1_\ZZ)-1=\rho(F)+1= \rho(H_e)=\rho(\cyc(H_e))+|H_e-\cyc(H_e)|=\rho(\cyc(H_e))+|F-X_e|,$$ where the last equality follows from Lemma~\ref{lm:FH}.
Thus by definition, $X_e\in\Upsilon(F)$. 
 
\underline{$\cyc(H_e)\in \bar{\Upsilon}(F)$:} Assume for a contradiction that $\cyc(H_e)$ is not inclusion-maximal in $\Upsilon(F)$, and that $K\in \Upsilon(F)$ is a cyclic flat with $\cyc(H_e)\subseteq K$ and $$\rho(1_\ZZ)-1=\rho(K)+|F-K|.$$ Then $H_e-K=F-K$, because by Lemma~\ref{lm:FH} $H_e-F\subseteq \cyc H_e\subseteq K$. Moreover, as $H_e$ is a flat, $$K\not\in\{\cyc H_e\}=\ZZ(H_e),$$ so $\rho(K)+|H_e-K|>\rho(H_e).$ This yields the chain of inequalities $$\rho(1_\ZZ)-1=\rho(K)+|F-K|=\rho(K)+|H-K|>\rho(H_e)=\rho(F)+1,$$ which contradicts the assumption $\rho(1_\ZZ)=\rho(F)+2$. Thus $\cyc(H_e)$ is inclusion-maximal in $\Upsilon(F)$.

\underline{$\bar{\Upsilon}(F)\subseteq\{\cyc(H_e) : e\in 1_\ZZ- F^\FF\}$:} Let $K\in \bar{\Upsilon}(F)$. If $K\subseteq F$, then as $K$ is cyclic we would have $$K\subseteq \cyc (F)\subseteq \cyc(\cl(F\cup e))\in\Upsilon(F)$$ for any $e\in 1_\ZZ- F^\FF$. This contradicts the maximality of $K$. 

We thus have $K-F\neq \emptyset$, or in other words $K\cap F\subsetneq K$, so $\rho(K\cap F)<\rho(F)$. It follows that $$\rho(1_ZZ)-2=\rho(F)\leq \rho(K\cap F)+|F-K|\leq \rho(K)+|F-K|-1=\rho(1_\ZZ) -2.$$ In particular we have equality $\rho(F)= \rho(K\cap F)+|F-K|$, meaning that $\cyc(F)\subseteq K\cap F$. 

As $K$ is cyclic, we have $|K-F|\geq 2$. On the other hand we have $$\rho(F)+1=\rho(1_\ZZ)-1=\rho(K)+|F-K|\geq \rho(K)+\rho(F)-\rho(F\cap K),$$ so $\rho(K)=\rho(K\cap F)+1$. It follows that for any $e\in K-F$, $F\cup e\not\in\FF$ and $K\subseteq \cl(F\cup e)=H_e$. Since $K$ is cyclic, we also get $K\subseteq\cyc(H_e)\in \Upsilon(F)$. But $K$ was assumed to be maximal in $\Upsilon(F)$, which shows that $K\subseteq\cyc(H_e)$.
\end{proof}

\begin{theorem}\label{thm:u2n}
Let $M=(E, \rho)$ be a matroid and let $F\subseteq 1_{\ZZ}$ be a flat with $\rho(F)=\rho(1_\ZZ)-2$. Then $$|(F, 1_\ZZ)_\FF| = |1_\ZZ - F| - \sum_{Z\in \bar{\Upsilon}(F)}(|Z-F|-1).$$
\end{theorem}
\begin{proof}
There is a natural surjective map $1_\ZZ-F\to (F, 1_\ZZ)_\FF$ given by $e\mapsto\cl(F\cup e)$. This map is injective on $F^\FF-F$, because for $e\in F^\FF-F$ we have $F\cup e=\cl(F\cup e)$.

%Let $(F, 1_\ZZ)_\FF = \{F_1,\ldots,F_m\}$, and observe that
%$$
%\{F_1 - F, \ldots, F_m - F \}\hbox{ is a partition of } 1_\ZZ-F.
%$$
%Thus every $F_i$ can be written as $\cl(F \cup e)$, for some $e\in 1_\ZZ-F$. 

So $|(F, 1_\ZZ)_\FF|-|F^\FF-F|$ is the number of sets $H_e=\cl(F\cup e)$ where $e\in 1_\ZZ-F^\FF$, and each such set corresponds to $|H_e-F|$ elements of $1_\ZZ-F^\FF$.
By Proposition~\ref{lm:upsilon}, the cyclic operator is a bijection from this collection to $\bar{\Upsilon}(F)$, so $$|(F, 1_\ZZ)_\FF|-|F^\FF-F|=|\{H_e: e\in 1_\ZZ-F^\FF\}|=|\bar{\Upsilon}(F)|.$$ By Lemma~\ref{lm:FH}, we also have $|\cyc(H_e)-F|=|H_e-F|$ for $e\in 1_\ZZ-F^\FF$, so $$\sum_{Z\in \bar{\Upsilon}(F)}|Z-F|=| 1_\ZZ-F^\FF|,$$ and hence $$\sum_{Z\in \bar{\Upsilon}(F)}|(Z-F)-1|=| 1_\ZZ-F^\FF\|-|\bar{\Upsilon}(F)|.$$

 Combining this, we get $$|(F, 1_\ZZ)_\FF|=|F^\FF-F|+|\bar{\Upsilon}(F)|=| 1_\ZZ-F|-\sum_{Z\in \bar{\Upsilon}(F)}|(Z-F)-1|.$$

\end{proof}

The dual version of Theorem~\ref{thm:u2n} identifies $U_{n}^{n-2}$ minors of $M$ via a reconstruction of $\UU(M)$. In particular, this yields two different ways to characterise $U_4^2$-avoiding matroids, by either reconstructing the upper part of $\FF(M)$, or by reconstructing the lower part of $\UU(M)$. 

\begin{corollary}~\label{cor:binary}
Let $M$ be a matroid. The following three conditions are equivalent.
\begin{enumerate}[(i)]
\item $M$ is binary.
\item For every flat $F\subseteq 1_\ZZ$ with $\rho(F)=\rho(1_Z)-2$ it holds that 
$$|1_\ZZ-F|-\sum_{X}(|X-F|-1)<4,$$ 
where the sum is taken over all $X\in \ZZ(M)$ such that $\rho(1_{\ZZ(M)})-1=\rho(X)+|F-X|$.
\item For every cyclic set $U\supseteq 0_\ZZ$ with $\eta(U)=|0_\ZZ|+2$ it holds that 
$$|U-0_\ZZ|-\sum_{X}(|U-X|-1)<4,$$
where the sum is taken over all $X\in \ZZ(M)$ such that $|0_\ZZ|+1+\rho(X)=|X\cap U|$. 
\end{enumerate}
\end{corollary}

 % F(M) from Z(M) + Char U_2n avoiding matroids

%Binary part-------------------
%Input "zurich" + "castle2"
%Contains :
%		Rank/nullity description  of Z(M) and config 
%		Atomicity of binary matroids.

%------------------------------------------------------------

\section{Covering Relations in $\ZZ(M)$ and Atomicity}

For the rest of the paper, we focus our study on binary matroids and their lattice of cyclic flats. The first consequence of restricting to binary matroids deals with the covering relations in the lattice of cyclic flats. By Theorem \ref{tutte}, a matroid is representable over $\F_2$ if and only if it avoids $U_4^2$ as a minor. In this case,  Corollary~\ref{corr:hasse} tells us that for $X \lessdot_{\ZZ} Y$ we cannot simultaneously have $\rho(Y)-\rho(X)>1$ and $\eta(Y)-\eta(X)>1$. 
On the other hand, we know by Theorem 3.2 in \cite{bonin08} or by direct calculation, that we always have $\rho(Y)-\rho(X)\geq 1$ and $\eta(Y)-\eta(X)\geq 1$. Thus, if $M$ is representable over $\F_2$, then every edge $X\lessdot_\ZZ Y$ in the Hasse diagram of $\ZZ(M)$ satisfies exactly one of the following:

\begin{enumerate}[(i)]
\item $\rho(Y)-\rho(X)=l>1$. We call such an edge a \emph{rank edge}, and label it $\rho=l$. Such an edge corresponds to a $U_{l+1}^l$ minor in $M$.
\item $\eta(Y)-\eta(X)=l>1$. We call such an edge a \emph{nullity edge}, and label it $\eta=l$. Such an edge corresponds to a $U_{l+1}^1$ minor in $M$.
\item $\rho(Y)-\rho(X)=1$ and $\eta(Y)-\eta(X)=1$. We call such an edge an \emph{elementary edge}. Such an edge corresponds to a $U_2^1$ minor in $M$.
\end{enumerate}

We illustrate this phenomenon in an example. 

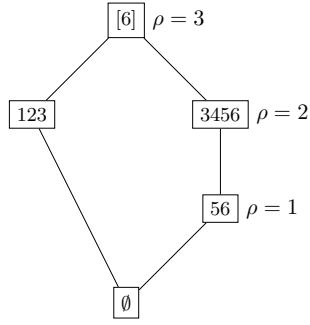
\begin{figure}
\centering
\resizebox{0.33\textwidth}{!}{%
\begin{tikzpicture}

	\node[shape=rectangle,draw=black] (00) at (0,0) {\small$\emptyset$};  
    \node[shape=rectangle,draw=black, label=right:{$\rho = 1$}] (12) at (1.5,1.5) {\small\(56\)}; 
    \node[shape=rectangle,draw=black] (21) at (-1.5,3) {\small\(123\)}; 
    \node[shape=rectangle,draw=black, label=right:{$\rho = 2$}] (22) at (1.5,3) {\small \(3456\)}; 
    \node[shape=rectangle,draw=black,, label=right:{$\rho = 3$}] (31) at (0,4.5) {\small \([6] \)};

     \path [-] (00) edge (12);
    \path [-] (00) edge (21);
    \path [-] (12) edge (22);
    \path [-] (21) edge (31);
    \path [-] (22) edge (31);

\end{tikzpicture}
}
\caption{Lattice of cyclic flats of the matroid from Example \ref{ex:ranknulledges}.}
\label{fig:LCF_ExR}
\end{figure}

\begin{example}
\label{ex:ranknulledges}
Let $M=([6],\rho)$ be the binary matroid generated by the matrix
\[
G=
\begin{pmatrix}
1&0&1&0&1&1\\
0&1&1&0&1&1\\
0&0&0&1&1&1\\
\end{pmatrix}.
\]
The Hasse diagram of its lattice of cyclic flats is displayed in Figure \ref{fig:LCF_ExR}. We have that $\rho(\{1,2,3\})=2$ and $\eta(\{1,2,3\})=1$. Therefore, in $\ZZ(M)$, the covering relation $\emptyset \lessdot \{1,2,3\}$ is a rank edge. On the other hand, since $\eta([6])=3$, the covering relation $\{1,2,3\} \lessdot [6]$ is a nullity edge. Finally, one the right-hand side of the Hasse diagram, every covering relation in the chain $\emptyset \lessdot \{5,6\} \lessdot \{3,4,5,6\} \lessdot [6] $ is an elementary edge. 
\end{example}

As we will see in the next sections of this paper, the dual relation of being a rank or a nullity edge plays a crucial role in understanding the lattice of cyclic flats. This relation also affects the possible parameters of a matroid and in particular its minimum distance. The study of the connection between the nullity edges and the minimum distance of a matroid is the topic of Section \ref{sec:recursive_structure}. We give here a first glimpse of this connection.  

\begin{lemma}
\label{lemma:d3_nulledge}
If $M=(E,\rho)$ is a binary matroid with minimum distance $d \geq 3$ then every edge $Z \lessdot E$ is a nullity edge.  
\end{lemma}

\begin{proof}
Since the minimum distance is greater than $1$, this implies that $M$ contains no isthmuses and $1_{\ZZ} = E$. Now the minimum distance satisfies the relation
\[
d = \eta(E) + 1 - \max \{ \eta(Z) : Z \in \ZZ(M) \st Z \lessdot E \}.
\]
Therefore, if $Z \in \ZZ(M)$ is such that $Z \lessdot E$, we have that $2 \leq \eta(E) - \eta(Z)$ and the edge $Z \lessdot E$ is a nullity edge. 
\end{proof}

Avoiding a $U_{4}^{2}$ minor is not the only characterization of binary matroids. In fact, it was proven in \cite{oxley11} that a binary matroid can be characterised by the relation between its circuits and the circuits of its dual matroid or directly by the symmetric difference of its circuits. 

\begin{theorem}[\cite{oxley11}]
\label{thm:char_binary_circuits}
Let $M$ be a matroid. The following are equivalent
\begin{enumerate}
\item $M$ is binary.
\item Let $C$ and $C^{\ast}$ be a circuit and a co-circuit respectively. Then $|C \cap C^{\ast}| $ is even.
\item Let $C_{1}, \ldots , C_{k}$ be circuits. Then the symmetric difference $C_{1} \bigtriangleup \ldots \bigtriangleup C_{k}$ is a disjoint union of circuits. 
\end{enumerate}
\end{theorem}

This naturally leads us to consider the relation between circuits and cyclic flats of a binary matroid. For this, we restrict our study to simple matroids. The reason is that parallel elements of a binary matroid are just repeated elements. We begin by an immediate consequence of the rank edges on the atoms of $\ZZ(M)$.

%Atomicity of \ZZ(M)-------
\begin{lemma}
\label{lemma:atom_iff_eta1}
Let $M=(E,\rho)$ be a simple binary matroid. Then $Z$ is an atom of $\ZZ(M)$ if and only if $\eta(Z)=1$. 
\end{lemma}

\begin{proof}
Since $M$ is simple, it guarantees that $\emptyset = 0_{\ZZ}$. Furthermore, it also means that, for all cyclic flats $Z \neq \emptyset$, we have $\rho(Z)>1$. Hence, every atom $Z_{at}$ will have a rank edge, \emph{i.e.}, $\eta(Z_{at}) = 1$. 
\end{proof}

The next two lemmas link atoms in $\ZZ(M)$ to certain minimal circuits in $M$. 

\begin{lemma}
\label{lemma:circuit_atom}
Let $M=(E,\rho)$ be a simple binary matroid. Let $C$ be a circuit of $M$. Then $\cl(C)$ is an atom of $\ZZ(M)$ if and only if $\cl(C)=C$. 
\end{lemma}

\begin{proof}
By lemma \ref{lemma:atom_iff_eta1}, $\cl(C)$ is an atom of $\ZZ(M)$ if and only if $\eta(\cl(C))=1$. But since $C$ is a circuit, we have $\eta(C)=1$. Now $\eta(\cl(C)) = 1 + |\cl(C)-C|$. Hence $\cl(C)$ is an atom if and only if $\cl(C) = C$.
\end{proof}

\begin{lemma}
\label{lemma:circuit_closure}
Let $M=(E,\rho)$ be a simple binary matroid that contains no isthmuses. Let $e \in E$ and $C$ be a circuit of minimal length containing $e$. Then $\cl(C)=C$.
\end{lemma}

\begin{proof}
First, the existence of a circuit containing $e$ is guaranteed by the fact that $M$ contains no isthmuses. Let $C$ be a minimum length circuit containing $e$. 

Now, consider a binary representation $\{x_{f} \}_{f \in M}$ of $M$. We can express $x_{e}$ by a linear combination of elements $\{x_{f} : f \in C\setminus \{e\} \}$. Since $C$ is a binary circuit, we will need all elements in $C \setminus \{e\}$ with coefficients equal to 1. Hence
\[
x_{e}=\sum\limits_{f \in C\setminus \{ e\} } x_{f}.
\]
Assume for a contradiction that there exists $e' \in \cl(C) - C$. Then $x_{e'} = \sum\limits_{f \in D \subseteq C}x_{f}$. Since $M$ is binary and simple, we have $2 \leq |D| < |C|$. 

If $e \in D$, then we have found a circuit smaller than $C$ containing $e$, which is a contradiction to the minimality of $C$. 

If $e \notin D$, then
\[
x_{e}= \sum\limits_{f \in D}x_{f} + \sum\limits_{f' \in (C\setminus \{e \} ) \setminus D} x_{f'} = x_{e'} + \sum\limits_{f' \in (C\setminus \{e \} ) \setminus D} x_{f'}.
\]
Thus, $e$ is in the circuit $\{e \} \cup \{e' \} \cup ((C\setminus \{e \}) \setminus D)$ with cardinality $| \{e'\} \cup (C \setminus D)| < |C|$ by the fact that $|D|\geq 2$. Again, this is a contradiction to the minimality of $C$. Hence $\cl(C)=C$.
\end{proof}

By combining Lemma \ref{lemma:circuit_atom} and \ref{lemma:circuit_closure}, we obtain the following result. 

\begin{lemma}
\label{lemma:element_in_atom}
Let $M=(E,\rho)$ be a simple binary matroid that contains no isthmuses. Then every element $e \in E$ is contained in an atom. 
\end{lemma}

We have enough results to prove the main result of this section, which states that the lattice of cyclic flats of a simple binary matroid with no isthmuses is atomic. 

\begin{theorem}
\label{thm:lcf_atomic}
Let $M=(E,\rho)$ be a simple binary matroid that contains no isthmuses. Then the lattice of cyclic flats $\ZZ(M)$ is atomic.
\end{theorem}

\begin{proof}
By Lemma \ref{lemma:element_in_atom}, for every $e \in E$, there exists an atom $Z_{at}^{e} \in \ZZ(M)$ with $e \in Z_{at}^{e}$. Thus, 
\[
\bigvee\limits_{e \in M} Z_{at}^{e} \supseteq \bigcup\limits_{e \in M}Z_{at}^{e} = E.
\]
For a cyclic flat $Y \in \ZZ(M)$, we can restrict the matroid to $M|Y=(Y,\rho)$. Since $\ZZ(M|Y) = \{ Z: Z \subseteq Y, Z \in \ZZ(M) \}$ by Corollary \ref{corr:interval} and $M|Y$ contains no co-loops, we are back to the previous case. Hence
\[
Y = \bigvee\limits_{e \in Y} Z_{at}^{e}
\]
and this proves that $\ZZ(M)$ is atomic. 
\end{proof}

Indeed, we have proven a slightly stronger property than atomicity. Namely, any element in $\ZZ(M)$ is equal not only to the join, but also to the union of all the atoms that it contains. As we can see in Example \ref{ex:ranknulledges}, it is crucial that the matroid is simple for the lattice of cyclic flats to be atomic. As a corollary, we obtain that the lattice of cyclic flats of a binary non-degenerate matroid is coatomic if the minimum distance is greater than $2$. 

\begin{corollary}
\label{cor:coatomic}
Let $M=(E,\rho)$ be a binary non-degenerate matroid. If the minimum distance $d \geq 3$, then $\ZZ(M)$ is coatomic. 
\end{corollary}

\begin{proof}
$M$ being non-degenerate implies that $M^{\ast}$ is also non-degenerate. Let $Z^{\ast}$ be an atom of $\ZZ(M^{\ast})$. By dual property, we have that $E-Z^{\ast}$ is a coatom of $\ZZ(M)$. Now Lemma \ref{lemma:d3_nulledge} implies that $\rho^{\ast}(Z^{\ast}) = \eta(E) - \eta(E-Z^{\ast}) \geq 2$. Hence $M^{\ast}$ contains no parallel elements and Theorem \ref{thm:lcf_atomic} implies that $\ZZ(M^{\ast})$ is atomic. 
\end{proof}

Finally, by combining the previous results, we obtain a relation between the atoms and the coatoms of $\ZZ(M)$.

\begin{lemma}
\label{lemma:atom_coat_relation}
Let $\rmatroid$ be a binary simple matroid with no isthmuses and $d \geq 3$. Then for every atom $Z_{a}$ and coatom $Z^{c}$ of $\ZZ(M)$ we have that $|Z_{a} \setminus Z^{c}|$ is even. 
\end{lemma}

\begin{proof}
Every atom of $\ZZ(M)$ is a circuit of $M$ by Lemma \ref{lemma:atom_iff_eta1}. Since $d\geq 3$, $M^{\ast}$ is simple and every coatom of $\ZZ(M)$ is the complement of a cocircuit. Therefore, by Theorem \ref{thm:char_binary_circuits}, we have that $|Z_{a} \cap (E-Z^{c})| = |Z_{a} \setminus Z^{c}|$ is even. 
\end{proof}

 % Rank / nullity description of Z(M) + Atomicity
%Input "height"
%Contains : Height 3 sublattices with rank 2 + Height 3 sublattices with d>2.

\section{Matroids with Lattices of Cyclic Flats of Height 3}
\label{sec:height3}

In this section, we study binary matroids when their lattice of cyclic flats has height 3. Under this assumption, every atom of $\ZZ(M)$ is also a coatom, which makes the structure of $\ZZ(M)$ very rigid. First, we focus on matroids of rank 2 and derive formulas that relate the nullity of the ground set, the number of atoms and the nullity of these atoms. Although technical, these formulas will be very useful in the next section when we study recursive structures in the lattice of cyclic flats. 

Secondly, we extend the study of height-three lattices to binary matroids with arbitrary rank. It turns out that only a few binary simple matroids can have a lattice of cyclic flats of height 3. In this part, we prove the non-existence of simple matroids with lattice of cyclic flats of height 3 depending on the size and rank and give the complete classification of these matroids when $\eta(E)$ is greater than or equal to $3$.

\subsection{Matroids of Rank 2 or Nullity 2}

We start by considering matroids of rank 2. For binary matroids, $U_{3}^{2}$ is the unique simple matroid of rank 2 that contains no isthmuses and has a lattice of cyclic flats of height 2. If the nullity of $M$ is larger than 1, it implies that some elements are parallel and thus the lattice of cyclic flats has height 3. Using this fact, we can express the nullity of $M$ depending on the nullity of the atoms and the number of atoms. 

\begin{lemma}
\label{lemma:rank2_null_relations}
Let $\rmatroid$ be a binary non-degenerate $(n,k,d)$-matroid with rank $k=2$. Let $\Upsilon_{\emptyset}$ be the set of atoms of $\ZZ(M)$. Then, we have the following relations:
\begin{itemize}
\item If $|\Upsilon_{\emptyset}|=0$, then $\eta(E)=1$. 
\item If $|\Upsilon_{\emptyset}|=1$, then $\eta(E)=\eta(Z) + 1$ with $\{Z\} = \Upsilon_{\emptyset}$. 
\item If $|\Upsilon_{\emptyset}|=2$ and $E=\bigcup\limits_{Z \in \Upsilon_{\emptyset}} Z$, then $\eta(E) = \sum\limits_{Z \in \Upsilon_{\emptyset}} \eta(Z)$.
\item If $|\Upsilon_{\emptyset}|=2$ and $E-\bigcup\limits_{Z \in \Upsilon_{\emptyset}} Z \neq \emptyset$, then $\eta(E) = 1 + \sum\limits_{Z \in \Upsilon_{\emptyset}} \eta(Z)$.
\item if $|\Upsilon_{\emptyset}|=3$, then $\eta(E) = 1 + \sum\limits_{Z \in \Upsilon_{\emptyset}} \eta(Z)$.
\end{itemize}
\end{lemma}

\begin{proof}
Let $G$ be the matrix associated to the matroid $M$. Since $G$ is a binary matrix of rank 2, there are only three possible choices for the columns of $G$, namely the vectors ${1 \choose 0}$, ${0 \choose 1}$ and ${1 \choose 1}$. Then, the size of $E$ can be counted as $|E| = \# {1 \choose 0} + \# {0 \choose 1} + \# {1 \choose 1}$. Furthermore, every time one vector is repeated, it will create a cyclic flat of rank 1 in $\ZZ(M)$. Thus, if $\Upsilon_{\emptyset}$ is the set of these cyclic flats, we have $|E|=\sum\limits_{Z \in \Upsilon} |Z| + |E- \bigcup\limits_{Z \in \Upsilon} Z|$. Since we also know their rank, we can transform the previous equation into an equation on the nullity. Now splitting this equation depending on the value of $|\Upsilon_{\emptyset}|$ will give the result (notice that since we assume no isthmuses, $|\Upsilon_{\emptyset}|=1$ forces the two other vectors to appear in $G$). 
\end{proof}

More interestingly, the previous formulas can be generalized as local relations on arbitrary binary matroids. Indeed, if two cyclic flats have a rank difference of 2, we can use contraction and deletion to obtain a rank $2$ matroid and apply Lemma \ref{lemma:rank2_null_relations}. Furthermore, by minor properties, these relations can be directly expressed in $M$ instead of in the minor obtained from $M$. 

\begin{lemma}
\label{lemma:rankdiff2_null_relations}
Let $\rmatroid$ be a binary matroid and let $Z_{1}, Z_{2} \in \ZZ(M)$ such that $Z_{1} \subset Z_{2}$ and $\rho(Z_{2})-\rho(Z_{1})=2$. Define $\Upsilon = \{ Z : Z \in \ZZ(M), \,  Z_{1} \lessdot Z \lessdot Z_{2} \}$. Then, we have the following relations. 
\begin{itemize}
\item If $|\Upsilon|=0$, then $\eta(Z_{2})=\eta(Z_{1})+1$. 
\item If $|\Upsilon|=1$, then $\eta(Z_{2})=\eta(Z) + 1$ with $\{Z\} = \Upsilon$.
\item If $|\Upsilon|=2$ and $Z_{2}=\bigcup\limits_{Z \in \Upsilon} Z$, then $\eta(Z_{2}) = \sum\limits_{Z \in \Upsilon} \eta(Z) - \eta(Z_{1})$.
\item If $|\Upsilon|=2$ and $Z_{2}-\bigcup\limits_{Z \in \Upsilon} Z \neq \emptyset$, then $\eta(Z_{2}) = 1 + \sum\limits_{Z \in \Upsilon} \eta(Z) - \eta(Z_{1})$.
\item If $|\Upsilon|=3$, then $\eta(Z_{2})=1+ \sum\limits_{Z \in \Upsilon} \eta(Z) -2\eta(Z_{1}) $.
\end{itemize}
\end{lemma}

\begin{proof}
The minor $M|Z_{2}/Z_{1}$ is a binary non-degenerate matroid. Hence we can apply Lemma \ref{lemma:rank2_null_relations} with ground set $Z_{2}-Z_{1}$ and nullity function $\eta_{M|Z_{2}/Z_{1}}=\eta_{M/Z_{1}}$. Now if $A \subset Z_{2}-Z_{1}$, then $\eta_{M|Z_{2}/Z_{1}}(A)=\eta(A \cup Z_{1}) - \eta(Z_{1})$. Using this with Lemma \ref{lemma:rank2_null_relations} will give the result.
\end{proof}

The next lemma is the dual version of the previous lemma. 

\begin{lemma}
\label{lemma:null2_rank_relations}
Let $\rmatroid$ be a binary matroid and let $Z_{1}, Z_{2} \in \ZZ(M)$ be such that $Z_{1} \subset Z_{2}$ and $\eta(Z_{2})-\eta(Z_{1})=2$. Define $\Upsilon = \{ Z : Z \in \ZZ(M), \, Z_{1} \lessdot Z \lessdot Z_{2} \}$. Then, we have the following relations.
\begin{itemize}
\item If $|\Upsilon|=0$, then $\rho(Z_{2})=\rho(Z_{1})+1$.
\item If $|\Upsilon|=1$, then $\rho(Z)=\rho(Z_{1})+1$ with $\{Z\} = \Upsilon$.
\item If $|\Upsilon|=2$ and $Z_{1}=\bigcap\limits_{Z \in \Upsilon} Z$, then $\sum\limits_{Z \in \Upsilon} \rho(Z) = \rho(Z_{1}) + \rho(Z_{2})$.
\item If $|\Upsilon|=2$ and $\left(\bigcap\limits_{Z \in \Upsilon} Z\right) - Z_{1} \neq \emptyset$, then $\sum\limits_{Z \in \Upsilon} \rho(Z) = 1 + \rho(Z_{1}) + \rho(Z_{2})$.
\item If $|\Upsilon|=3$, then $\sum\limits_{Z \in \Upsilon} \rho(Z) = 1 + \rho(Z_{1}) + \rho(Z_{2})$. 
\end{itemize}
\end{lemma}

\subsection{Matroids of Arbitrary Rank}

In this part, we relax the condition on the rank while still forcing the lattice of cyclic flats to have height $3$. We will see that, in fact, only a few simple binary matroids satisfy this condition and we will completely characterize them. We first treat the case when the matroids contain parallel elements. 

\begin{proposition}
Let $\rmatroid$ be a binary non-degenerate $(n,k,d)$-matroid with $k\geq 3$. If $\ZZ(M)$ has height 3 and $M$ contains parallel elements, then $d=2$. 
\end{proposition}

\begin{proof}
Let $e \in E$ be one of the parallel elements. Since $\ZZ(M)$ has height $3$, the lattice of cyclic flats contains the chain $\emptyset \lessdot \cl(e) \lessdot E$. Now $\rho(\cl(e))=\rho(e)=1$. Then we have that $\cl(e) \lessdot E$ is a rank edge, implying that $\eta(E)=\eta(\cl(e))+1$. Hence $d=\eta(E) + 1 - \max \{ \eta(Z) : Z \in \ZZ(M)-E \} \leq \eta(E) + 1 - \eta(\cl(e)) = 2$, and since $M$ is non-degenerate, we have $d=2$. 
\end{proof}

Thus, it is always possible to increase the nullity of $M$ by adding parallel elements. However, if $\ZZ(M)$ has height 3, then the minimum distance is always equal to $2$. We focus now on simple matroids and start by an upper bound on the intersection between two atoms. 

\begin{lemma}
\label{lemma:atoms_intersection}
Let $\rmatroid$ be a binary simple $(n,k,d)$-matroid with no isthmuses and $k \geq 3$. If $\ZZ(M)$ has height 3 and for every atom $Z \in \ZZ(M)$ we have $\rho(Z)=k-1$, then 
\[
|Z_{1} \cap Z_{2}| \leq \frac{k}{2}, \text{ for all } Z_{1}, Z_{2} \in \ZZ(M)-E.
\]
\end{lemma}

\begin{proof}
Let $Z_{1}, Z_{2} \in \ZZ(M)-E$. If one of them is the empty set, then the result is trivial. Assume now that neither of them are empty. Since $\ZZ(M)$ has height 3, $Z_{1}$ and $Z_{2}$ must be atoms of $\ZZ(M)$ with parameters $\rho(Z_{i})=k-1$, $\eta(Z_{i})=1$ and $|Z_{i}|=k$ for $1 \leq i \leq 2$. Furthermore, $Z_{1}$ and $Z_{2}$ are also circuits in $M$. By Theorem \ref{thm:char_binary_circuits}, $Z_{1} \bigtriangleup Z_{2}$ is a disjoint union of circuits. Using the fact that $Z_{1} \cup Z_{2} = Z_{1} \bigtriangleup Z_{2} \uplus Z_{1} \cap Z_{2}$, we have that
\[
|Z_{1} \bigtriangleup Z_{2}| = |Z_{1} \cup Z_{2}| - |Z_{1} \cap Z_{2}| = |Z_{1}| + |Z_{2}| - 2 |Z_{1} \cap Z_{2}| = 2k - 2 |Z_{1} \cap Z_{2}|.
\]

Now the smallest size of a circuit in $M$ is $k$ since otherwise a circuit of size less than $k$ will yield a cyclic flat of $\ZZ(M)$ of rank less than $k-1$, and thus contradict our assumptions. This implies that $|Z_{1} \bigtriangleup Z_{2}| \geq k$ and hence $|Z_{1} \cap Z_{2}| \leq \frac{k}{2}$. 
\end{proof}

We now prove one of the main results of this section. The next proposition states the non-existence of simple matroids with lattice of cyclic flats of height 3 and large rank and nullity. In fact, as soon as the rank is larger than or equal to $5$, the only possible such matroids have nullity $2$ and thus need to satisfy the relations in Lemma \ref{lemma:null2_rank_relations}. 

\begin{proposition}
\label{prop:3l_cf_k5}
Let $\rmatroid$ be a binary simple $(n,k,d)$-matroid with no isthmuses. If $\ZZ(M)$ has height 3 and $k\geq 5$, then $\eta(E) = 2$ and $d =2$. 
\end{proposition}

\begin{proof}
If there exist $Z \in \ZZ(M)-\emptyset$ with $\rho(Z) < k-1$, then $Z \lessdot E$ is a rank edge and $\eta(E)=\eta(Z)+1 = 2$.  

Assume now that every atom of $\ZZ(M)$ has rank $k-1$ and assume for a contradiction that $\eta(E)>2$. The goal is to use Theorem \ref{thm:u2n} and Corollary \ref{cor:binary} to obtain a contradiction on the fact that $M$ is binary. 

Let $Z_{a}$ be an atom of $\ZZ(M)$. Remember that $\rho(Z_{a}) = k-1$ and $|Z_{a}|=k$. Choose $F \subset Z_{a}$ with $|F|=k-2$. Since the smallest size of a circuit is $k$, we have that $F \in \FF(M)$ and is independent. Since $1_{\ZZ}=E$, we have
\begin{align*}
\Upsilon(F) &= \{ X \in \ZZ(M) : \rho(E)-1 = \rho(X) + |F-X| \} \\
	&= \{ X \in \ZZ(M)-E : k-1 = k-1 + |F-X| \} \\
	&= \{ X \in \ZZ(M)-E : F \subset X \}.
\end{align*}
By Lemma \ref{lemma:atoms_intersection}, if $X$ and $X'$ are both atoms of $\ZZ(M)$, then $|X \cap X'| \leq \frac{k}{2}$. Since $k\geq 5$, we have $k-2>\frac{k}{2}$. This implies that $Z_{a}$ is the unique cyclic flat different from $E$ that contains $F$. Hence $\Upsilon(F) = \{ Z_{a} \} $. By Theorem \ref{thm:u2n} and since $\eta(E)>2$, we have
\[
|(F,E)_{\FF}|=|E-F| - |Z_{a}-F|+1=n-k+1 \geq 4.
\]

Therefore, by Corollary \ref{cor:binary}, $M$ is not binary, which contradicts our assumption. Thus, we have $\eta(E) \leq 2$. Since there exist atoms with nullity equal to 1, we get $\eta(E)=2$. Finally, the minimum distance is given by $d = \eta(E) + 1 - \max \{ \eta(Z) : Z \in \ZZ(M) -E \} = 2$.
\end{proof}

The next proposition goes further by giving an upper bound on the nullity when the rank is equal to $4$. 

\begin{proposition}
\label{prop:3l_cf_n9}
Let $\rmatroid$ be a binary simple $(n,k,d)$-matroid with no isthmuses and $k=4$. If $\ZZ(M)$ has height 3 then $n \leq 8$ or equivalently, $\eta(E)\leq 4$. 
\end{proposition}

\begin{proof}
If there exist $Z \in \ZZ(M)-\emptyset$ with $\rho(Z) < k-1$, then $Z \lessdot E$ is a rank edge and $\eta(E)=\eta(Z)+1 = 2$. 

Assume now that every atom of $\ZZ(M)$ has rank $k-1$. Let $F \subset E$ such that $|F|=k-2$. Since the smallest size of a circuit is $k$, we have that $F \in \FF(M)$ and is independent. By the same argument as in the proof of Proposition \ref{prop:3l_cf_k5}, we get that 
\[
\Upsilon(F)=\{ X \in \ZZ(M)-E : F \subset X \}.
\] 
By Theorem \ref{thm:u2n} and since $k=4$, we have
\[
|(F,E)_{\FF}|=|E-F| - \sum\limits_{Z \in \Upsilon(F)}(|Z-F|-1) = n-2 - |\Upsilon(F)|.
\]
By Corollary \ref{cor:binary}, we need $n-2 - |\Upsilon(F)| \leq 3$ or equivalently $n \leq 5 + |\Upsilon(F)|$. On the other hand, we have $n \geq |F| + |\Upsilon(F)| \cdot 2$. By combining the two equations, we get
\[
2 + 2|\Upsilon(F)| \leq 5 + |\Upsilon(F)| \iff |\Upsilon(F)| \leq 3.
\]
Hence we obtain $n \leq 8$. 
\end{proof}

The previous propositions restrict the candidates for simple matroids with height-three lattices to $k \leq 4$ and $n\leq 8$. We pursue by studying the structure of the lattice of cyclic flats of simple matroids with feasible size and rank. In particular, we prove that matroids satisfying these conditions are unique up to isomorphism. 

\begin{proposition}
Up to isomorphism, there is a unique binary simple $(6,3,3)$-matroid with no isthmuses. 
\end{proposition}

\begin{proof}
Let $\rmatroid$ be a binary simple $(6,3,3)$-matroid with no isthmuses. We start by proving that $\ZZ(M)$ has a unique configuration by counting the number of atoms. 

Since $M$ is simple with $k=3$ and $\eta(E)=3$, $\ZZ(M)$ has height 3. Let $Z_{a}$ be an atom of $\ZZ(M)$. We have $\rho(Z_{a})= 2$, $\eta(Z_{a}) = 1$, and $|Z_{a}|=3$. As the size of an atom is odd, by Lemma \ref{lemma:atom_coat_relation} and Lemma \ref{lemma:atoms_intersection} for all atoms $Z_{1}$ and $Z_{2}$, we have $Z_{1} \cap Z_{2} = 1$. Since $| Z_{1} \cup Z_{2} | \leq 5$, but $|E|=6$ and there always exists an atom for every coordinate, the number of atoms is at least 3. 

Denote by $Z_{1}, Z_{2}$, and $Z_{3}$ the first three atoms. Notice that they have to intersect pairwise in a different element because $|E|=6$. We have $|Z_{1} \Delta Z_{2} \Delta Z_{3}| = 3 $. Then, by Theorem \ref{thm:char_binary_circuits} on the symmetric difference, these three elements form an extra atom. Hence, the number of atoms is at least $4$. Now the number of atoms is upper bound by ${6 \choose 2} / 3 = 5$ since every pair will define a unique atom with 3 elements. However we cannot have 5 atoms. Indeed, by the inclusion--exclusion principle, we would have
\[
\left\vert \bigcup\limits_{i=1}^{5} Z_{i} \right\vert = 5 |Z_{1}| - 10 |Z_{1} \cap Z_{2}| = 5
\]
because every triple has an empty intersection. But this is not possible since already $|Z_{1} \cup Z_{2} \cup Z_{3}| = |E|=6$. Hence, $\ZZ(M)$ has 4 atoms and has a unique configuration. 

Now suppose $E=\{a,b,c,d,e,f\}$ and $Z_{1}=\{a,b,c\}$. By the previous part, we have $|Z_{1} \cap Z_{i}|=1$ for all $i \in \{2,3,4\}$ with a different element for each intersection. So, let $a \in Z_{2}, b \in Z_{3}$ and $c \in Z_{4}$. We complete $Z_{2}$ with some of the remaining elements to get $Z_{2}=\{a,d,e\}$. Since $|Z_{2} \cap Z_{3}|=1$, we choose $d \in Z_{3}$. The only possible choice for the last element in $Z_{3}$ is therefore $f$ and we have $Z_{3}=\{b,d,f\}$. Finally $Z_{4}$ has a non-trivial intersection with $Z_{2}$ and $Z_{3}$ so it has to be $Z_{4}=\{c,e,f\}$. Hence we have uniquely reconstructed the lattice of cyclic flats up to a permutation of the groundset which implies that $M$ is unique up to isomorphism. 
\end{proof}

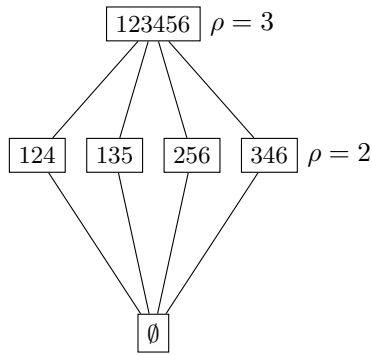
\begin{figure}[H]
\centering
\resizebox{0.4\textwidth}{!}{%
\begin{tikzpicture}

\node[shape=rectangle, draw=black] (00) at (0,0) {\small$\emptyset$};
 \node[shape=rectangle,draw=black] (10) at (-1.5,2.3) {\small\(124\)};
 \node[shape=rectangle,draw=black] (11) at (-0.5,2.3) {\small\(135\)};
 \node[shape=rectangle,draw=black] (12) at (0.5,2.3) {\small\(256\)};
 \node[shape=rectangle,draw=black, label=right:{$\rho = 2$}] (13) at (1.5,2.3) {\small\(346\)};
  \node[shape=rectangle,draw=black, label=right:{$\rho = 3$}] (20) at (0,4) {\small\(123456\)}; 

\path [-] (00) edge (10);
 \path [-] (00) edge (11);
 \path [- ] (00) edge (12);
 \path [- ] (00) edge (13);
\path [- ] (10) edge (20);
 \path [- ] (11) edge (20);
 \path [- ] (12) edge (20);
 \path [- ] (13) edge (20); 
    
\end{tikzpicture}
}
\caption{Lattice of cyclic flats of the binary $(6,3,3)$-matroid.}
\label{fig:633}
\end{figure}

Figure \ref{fig:633} displays the lattice of cyclic flats of the binary $(6,3,3)$-matroid with generator matrix 
\[
G=
\begin{pmatrix}
1&0&0&1&1&1\\
0&1&0&1&0&1\\
0&0&1&0&1&1\\
\end{pmatrix}.
\]

The second simple matroid of rank 3 that has a lattice of cyclic flats of height 3 is known as the simplex code $(7,3,4)$ where the generator matrix contains every possible column except the all-zero column. By definition, it is unique and has the maximum number of atoms which is ${7 \choose 2}/3 = 7$. Its dual is the $(7,4,3)$ Hamming code, which also has a lattice of cyclic flats of height 3 with 7 atoms. 

Finally, we study the last possible set of parameters. 
\begin{proposition}
There is a unique, up to isomorphism, binary simple $(8,4,4)$-matroid with no isthmuses. 
\end{proposition}

\begin{proof}
Let $\rmatroid$ be a binary simple $(8,4,4)$-matroid with no isthmuses. We again start by proving that $\ZZ(M)$ has a unique configuration by counting the number of atoms. 

Let $e \in E$. The number of atoms can be split into two groups: the atoms containing $e$ denoted by $A$ and the atoms not containing $e$ denoted by $B$. By Theorem \ref{thm:cf_formulas}, the number of atoms containing $e$ is greater than the number of atoms in $\ZZ(M/\{e\})$. Since $M/{e}$ is isomorphic to the $(7,3,4)$-matroid, we have that $|A| \geq 7$. Now $|B|$ is greater than the number of atoms in $\ZZ(M|(E-\{e\})$. The matroid $M|(E-\{e\})$ is isomorphic to the $(7,4,3)$-matroid which also contains $7$ atoms so $|B| \geq 7$. As the total number of atoms in an $(8,4,4)$-matroid cannot exceed ${8 \choose 3}/4 = 14$, $\ZZ(M)$ indeed contains $14$ atoms and has a unique configuration. Furthermore, $M$ is also the unique extension by one element of the Simplex $(7,3,4)$-matroid such that the contraction $M/e$ yields again the Simplex matroid and the deletion $M \setminus e$ yields the dual of this Simplex matroid, the $(7,4,3)$-matroid. 
\end{proof}

\begin{figure}[H]
\centering
\resizebox{0.9\textwidth}{!}{%
\begin{tikzpicture}

\node[shape=rectangle, draw=black] (00) at (0,0) {\small$\emptyset$};
\node[shape=rectangle,draw=black] (10) at (-6.5,4) {\small\(1235\)};
 \node[shape=rectangle,draw=black] (11) at (-5.5,4) {\small\(1248\)};
 \node[shape=rectangle,draw=black] (12) at (-4.5,4) {\small\(1267\)};
 \node[shape=rectangle,draw=black] (13) at (-3.5,4) {\small\(1347\)};
 \node[shape=rectangle,draw=black] (14) at (-2.5,4) {\small\(1368\)};
 \node[shape=rectangle,draw=black] (15) at (-1.5,4) {\small\(1456\)};
 \node[shape=rectangle,draw=black] (16) at (-0.5,4) {\small\(1578\)};
\node[shape=rectangle,draw=black] (17) at (0.5,4) {\small\(2346\)};
 \node[shape=rectangle,draw=black] (18) at (1.5,4) {\small\(2378\)};
 \node[shape=rectangle,draw=black] (19) at (2.5,4) {\small\(2457\)};
 \node[shape=rectangle,draw=black] (110) at (3.5,4) {\small\(2568\)};
 \node[shape=rectangle,draw=black] (111) at (4.5,4) {\small\(3458\)};
 \node[shape=rectangle,draw=black] (112) at (5.5,4) {\small\(3567\)};
 \node[shape=rectangle,draw=black, label=right:{$\rho = 3$}] (113) at (6.5,4) {\small\(4678\)};
 \node[shape=rectangle,draw=black, label=right:{$\rho = 4$}] (20) at (0,7) {\small\(12345678\)}; 

\path [-] (00) edge (10);
 \path [-] (00) edge (11);
 \path [- ] (00) edge (12);
 \path [- ] (00) edge (13);
 \path [- ] (00) edge (14);
 \path [- ] (00) edge (15);
 \path [- ] (00) edge (16);
 \path [- ] (00) edge (17);
 \path [- ] (00) edge (18);
 \path [- ] (00) edge (19);
 \path [- ] (00) edge (110);
 \path [- ] (00) edge (111);
 \path [- ] (00) edge (112);
 \path [- ] (00) edge (113);
\path [- ] (10) edge (20);
 \path [- ] (11) edge (20);
 \path [- ] (12) edge (20);
 \path [- ] (13) edge (20);
 \path [- ] (14) edge (20);
 \path [- ] (15) edge (20);
 \path [- ] (16) edge (20);
 \path [- ] (17) edge (20);
 \path [- ] (18) edge (20);
 \path [- ] (19) edge (20);
 \path [- ] (110) edge (20);
 \path [- ] (111) edge (20);
 \path [- ] (112) edge (20);
 \path [- ] (113) edge (20);
    
\end{tikzpicture}
}
\caption{Lattice of cyclic flats of the binary $(8,4,4)$-matroid.}
\label{fig:844}
\end{figure}
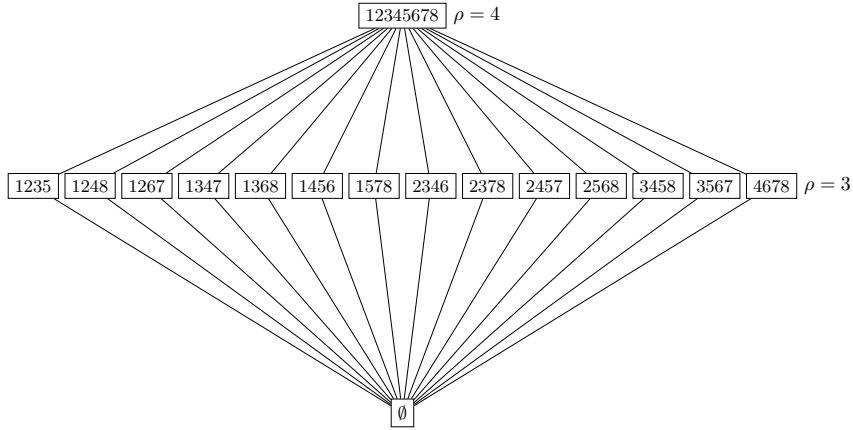

The binary code that satisfies these requirements can be obtained by considering the Reed--Muller code RM$(1,3)$ which indeed gives an $(8,4,4)$-matroid. The lattice of cyclic flats is displayed in Figure \ref{fig:844} with the generator matrix 
\[
G=
\begin{pmatrix}
1&0&0&0&1&0&1&1\\
0&1&0&0&1&1&0&1\\
0&0&1&0&1&1&1&0\\
0&0&0&1&0&1&1&1\\
\end{pmatrix}.
\]

We can now summarize the previous results. 

\begin{theorem}
\label{thm:classif_height3}
Let $\rmatroid$ be a binary simple $(n,k,d)$-matroid with no isthmuses. If $\ZZ(M)$ has height 3, then either $\eta(E)=2$ or $M$ is isomorphic to one of the matroids listed below:
\begin{itemize}
\item The $(6,3,3)$-matroid with 4 atoms.
\item The $(7,3,4)$-matroid with 7 atoms.
\item The $(7,4,3)$-matroid with 7 atoms.
\item The $(8,4,4)$-matroid with 14 atoms.
\end{itemize}
\end{theorem}

\begin{corollary}
\label{cor:height3_null4}
Let $\rmatroid$ be a binary simple matroid with no isthmuses. If $\ZZ(M)$ has height 3 then $\eta(E)\leq 4$ and $d\leq 4$. 
\end{corollary} % Matroids with height 3 Z(M)
%Input "griesmer"
%Contains : Recursive structure on coatom level. 

\section{Recursive Structure on Coatoms Level}
\label{sec:recursive_structure}

This section is devoted to understanding the consequences of the minimum distance to the shape of the top part of the lattice of cyclic flats. We already saw in Lemma \ref{lemma:d3_nulledge} that requiring the minimum distance to be greater than 2 forces all coatoms to have the same rank. Let us give this property a name for an arbitrary cyclic flat. 

\begin{definition}
A cyclic flat $Z$ of a binary matroid $M$ is \emph{blunt} if for all $Z' \in \ZZ(M)$ such that $Z' \lessdot Z$, we have $\rho(Z')=\rho(Z)-1$. In other words, every covering relation of a blunt cyclic flat is either a nullity edge or an elementary edge. 
\end{definition}

Thus, for a binary matroid $\rmatroid$, $E$ being blunt is a necessary condition to have $d \geq 3$. In this section, we study how the value of the minimum distance creates a recursive structure consisting of blunt cyclic flats in $\ZZ(M)$. In the second part, we state the equivalent notion of residual codes for binary matroids, which naturally leads to a version of the Griesmer bound for binary matroids. Finally, we discuss the relation between coatoms of $\ZZ(M)$ and codewords of the associated linear code. Before we begin, let us fix some notation.

\begin{notation}
We will usually denote a coatom by a superscripted $Z$ such as $Z^{1}$ or $Z^{c}$. If $\rmatroid$ is an $(n,k,d)$-matroid, then $Z^{d}$ denotes a cyclic flat with maximal nullity, \ie, we have $d=\eta(E) +1 - \eta(Z^{d})$. Notice that if $\rho(Z^{d})=k-1$, we have $|Z^{d}|=n-d$. 
\end{notation}

\subsection{Existence of Rank $k-2$ Cyclic Flats and Recursive Structure}

In order to get a recursive structure consisting of blunt cyclic flats, we are interested in the level below the coatoms level and in particular, in the cyclic flats of rank $k-2$. 

\begin{proposition}
\label{prop:ctrapos_double_nulledge}
Let $\rmatroid$ be a binary non-degenerate $(n,k,d)$-matroid and let $Z^{d} \in \ZZ(M)$ be a coatom with maximal nullity. If $Z^{d}$ is not blunt, then $d\leq 4$. 
\end{proposition}

\begin{proof}
First, notice that if $E$ is not blunt then $d=2$ by Lemma \ref{lemma:d3_nulledge}. Assume now that all coatoms have a rank equal to $k-1$. Since $Z^{d}$ is not blunt, there exists a $Z_{1} \in \ZZ(M)$ such that $Z_{1} \lessdot Z^{d}$ and $\rho(Z_{1}) < \rho(Z^{d})-1=k-2$.

Now the proof is a direct consequence of the classification of simple matroids with lattice of cyclic flats of height 3 in Section \ref{sec:height3}. Indeed, we will prove that $M/Z_{1}$ is simple with no isthmuses and its lattice of cyclic flats has height 3. 

Since $Z_{1} \lessdot Z^{d}$ is a rank edge, we have $\eta(Z^{d}) = \eta(Z_{1}) +1$. This implies that if $Z^{1}$ is a cyclic flat such that $Z_{1} \lessdot Z^{1}$, then $Z^{1}$ is a coatom with $\rho(Z^{1})=k-1$ because otherwise, there is a coatom $Z^{2}$ of $\ZZ(M)$ that covers $Z^{1}$ and $\eta(Z^{2}) \geq \eta(Z_{1}) + 2$, which contradicts the fact that $\eta(Z^{d})$ is maximal. Hence, by Theorem \ref{thm:cf_formulas}, $\ZZ(M/Z_{1})$ has height 3. 

The matroid $M/Z_{1}$ is also simple with no isthmuses as the contraction by a cyclic flat will not create any loops or isthmuses. Plus, for all coatoms $Z^{c} \in \ZZ(M)$ with $Z_{1} \lessdot Z^{c}$, we have $\rho(Z^{c}) - \rho(Z_{1}) \geq k-1 - (k-3) =2$, which means that there is no parallel elements in $M/Z_{1}$. Thus, $M/Z_{1}$ corresponds to the type of lattice studied in Section \ref{sec:height3}. 

Finally, to prove the proposition, we link the minimum distance $d$ to the minimum distance of $M/Z_{1}$. We have
\begin{align*}
d_{M/Z_{1}} & = \eta_{M/Z_{1}}(E-Z_{1}) + 1 - \max \{ \eta_{M/Z_{1}}(\tilde{Z}) : \tilde{Z} \in \ZZ(M/ Z_{1}) - (E-Z_{1}) \} \\
	& = \eta(E) - \eta(Z_{1}) + 1 - \max \{ \eta(Z-Z_{1}) : Z \in \ZZ(M)-E \nd Z_{1} \subseteq Z \} \\
	& = \eta(E) - \eta(Z_{1}) + 1 - \max \{ \eta(Z) : Z \in \ZZ(M)-E \nd Z_{1} \subseteq Z \} + \eta(Z_{1}) \\
	& = \eta(E) + 1 - \eta(Z^{d}) \\
	 & = d.
\end{align*}

Hence, by the classification obtained in Section \ref{sec:height3} and, in particular, by Corollary \ref{cor:height3_null4}, we have $d=d_{M/Z_{1}} \leq 4$.
\end{proof}

We now present two examples where no coatoms are blunt and the matroids have minimum distance $3$ and $4$ respectively.  

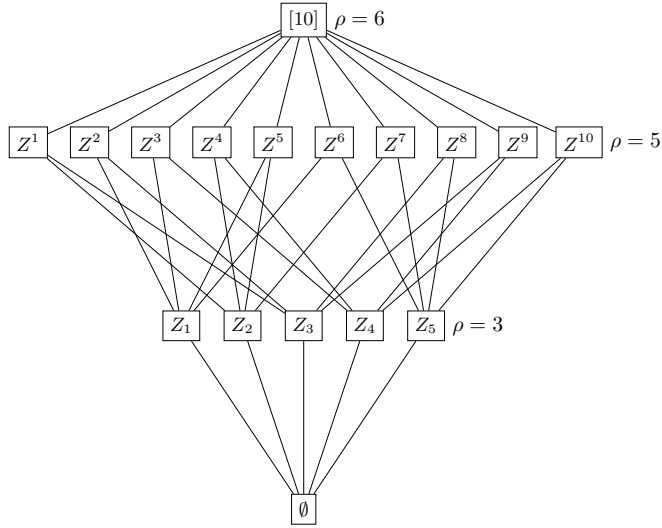
\begin{figure}
\centering
\resizebox{0.7\textwidth}{!}{%
\begin{tikzpicture}

\node[shape=rectangle,draw=black] (00) at (0,0) {\small$\emptyset$};

\node[shape=rectangle,draw=black] (10) at (-2,3) {\small $Z_{1}$ };
\node[shape=rectangle,draw=black] (11) at (-1,3) {\small $Z_{2}$};
\node[shape=rectangle,draw=black] (12) at (0,3) {\small $Z_{3}$};
\node[shape=rectangle,draw=black] (13) at (1,3) {\small $Z_{4}$};
\node[shape=rectangle,draw=black, label=right:{$\rho = 3$}] (14) at (2,3) {\small $Z_{5}$};

\node[shape=rectangle,draw=black] (20) at (-4.5,6) {\small $Z^{1}$};
\node[shape=rectangle,draw=black] (21) at (-3.5,6) {\small $Z^{2}$};
\node[shape=rectangle,draw=black] (22) at (-2.5,6) {\small $Z^{3}$};
\node[shape=rectangle,draw=black] (23) at (-1.5,6) {\small $Z^{4}$};
\node[shape=rectangle,draw=black] (24) at (-0.5,6) {\small $Z^{5}$};
\node[shape=rectangle,draw=black] (25) at (0.5,6) {\small $Z^{6}$};
\node[shape=rectangle,draw=black] (26) at (1.5,6) {\small $Z^{7}$};
\node[shape=rectangle,draw=black] (27) at (2.5,6) {\small $Z^{8}$};
\node[shape=rectangle,draw=black] (28) at (3.5,6) {\small $Z^{9}$};
\node[shape=rectangle,draw=black, label=right:{$\rho = 5$}] (29) at (4.5,6) {\small $Z^{10}$};
 
\node[shape=rectangle,draw=black, label=right:{$\rho = 6$}] (30) at (0,8) {\small $[10]$};
 
%paths----
\path [- ] (00) edge (10);
\path [- ] (00) edge (11);
\path [- ] (00) edge (12);
\path [- ] (00) edge (13);
\path [- ] (00) edge (14);

\path [-] (10) edge (21);
\path [-] (10) edge (22);
\path [-] (10) edge (24);
\path [-] (10) edge (25);

\path [-] (11) edge (20);
\path [-] (11) edge (23);
\path [-] (11) edge (24);
\path [-] (11) edge (26);

\path [-] (12) edge (20);
\path [-] (12) edge (21);
\path [-] (12) edge (27);
\path [-] (12) edge (28);

\path [-] (13) edge (22);
\path [-] (13) edge (23);
\path [-] (13) edge (28);
\path [-] (13) edge (29);

\path [-] (14) edge (25);
\path [-] (14) edge (26); 
\path [-] (14) edge (27);
\path [-] (14) edge (29); 

\path [-] (20) edge (30);
\path [-] (21) edge (30);
\path [-] (22) edge (30);
\path [-] (23) edge (30);
\path [-] (24) edge (30);
\path [-] (25) edge (30);
\path [-] (26) edge (30);
\path [-] (27) edge (30);
\path [-] (28) edge (30);
\path [-] (29) edge (30);
    
\end{tikzpicture}
}
\caption{Configuration of the lattice of cyclic flats of $M^{\ast}(K_{5})$.}
\label{fig:dual_K5}
\end{figure}

\begin{example}
Let $M$ be the binary $(10,6,3)$-matroid obtained by the dual of the complete graph $K_{5}$ and $G$ the following generator matrix of $M$ : 
\[
G=
\begin{pmatrix}
1&0&0&0&0&0&0&1&0&1\\
0&1&0&0&1&0&0&1&1&0 \\
0&0&1&0&1&0&0&0&1&0 \\
0&0&0&1&1&0&0&1&1&1 \\
0&0&0&0&0&1&0&1&1&0 \\
0&0&0&0&0&0&1&1&1&1 \\
\end{pmatrix}.
\]

We can see that there is no cyclic flat of rank $4$, since if $Z_{a}$ is an atom of $\ZZ(M)$ then the contracted matroid $M/Z_{a}$ is isomorphic to the $(6,3,3)$-matroid having a lattice of cyclic flats of height 3. The configuration of the lattice of cyclic flats is displayed in Figure \ref{fig:dual_K5}.
\end{example}

\begin{example}
The second example is the binary Reed--Muller code RM$(2,4)$ giving a $(16,11,4)$-matroid. Here, all atoms of $\ZZ(M)$ have rank $7$ and if $Z_{a}$ is one of them, then the contracted matroid $M/Z_{a}$ is isomorphic to the $(8,4,4)$-matroid, which has a lattice of cyclic flats of height 3. Thus, there are no cyclic flats of rank $9$. 
\end{example}

Now, we extend Proposition \ref{prop:ctrapos_double_nulledge} to coatoms with different size, $\ie$, bound the minimum distance when there is a coatom $Z^{c} \in \ZZ(M)$ which is not blunt. We emphasize that while coatoms $Z^{d} \in \ZZ(M)$ with maximal nullity always exist, coatoms with size less than $Z^{d}$ might not exist. For example, matroids coming from Simplex codes have only coatoms with maximal size. However, depending on the parameters $(n,k,d)$, we can use some techniques to guarantee the existence of smaller coatoms as demonstrated in Example \ref{ex:110405}.

We start by giving a lower bound on the number of coatoms of $\ZZ(M)$ covering a cyclic flat of rank $k-2$ when $d \geq 3$. 

\begin{lemma}
\label{lemma:upsilon2}
Let $\rmatroid$ be a binary non-degenerate $(n,k,d)$-matroid with $d \geq 3$, $Z_{1}$ a cyclic flat with $\rho(Z_{1})=k-2$ and $\Upsilon_{Z_{1}}$ the set of coatoms covering $Z_{1}$. Then we have $|\Upsilon_{Z_{1}}|\geq 2$. 
\end{lemma}

\begin{proof}
If $d\geq 3$, then Corollary \ref{cor:coatomic} states that $\ZZ(M)$ is co-atomic. Thus, there is at least two coatoms that cover a rank-$(k-2)$ cyclic flat. 
\end{proof}

We can now formulate an upper bound on $d$ when a coatom $Z^{c} \in \ZZ(M)$ is not blunt. The bound is expressed in terms of the gap between the nullity of $Z^{c}$ and the maximal nullity of a coatom, or equivalently, between their size difference. 

\begin{proposition}
\label{prop:bound_null_edge_d}
Let $\rmatroid$ be a binary non-degenerate $(n,k,d)$-matroid with $d \geq 3$ and $Z^{d} \in \ZZ(M)$ be such that $|Z^{d}|=n-d$. If there exists a coatom $Z^{c} \in \ZZ(M)$ not blunt with $|Z^{c}|<|Z^{d}|$, then
\[
d \leq 2 (\eta(Z^{d}) - \eta(Z^{c}) + 1) = 2 ( |Z^{d}|-|Z^{c}| +1 ).
\]
\end{proposition}

\begin{proof}
Since $Z^{c}$ is not blunt, there exists a cyclic flat $Z_{1} \lessdot Z^{c}$ with $\rho(Z_{1}) < k-2$. If there are no cyclic flats of rank $k-2$ that contain $Z_{1}$, then let $Z_{m}$ be the biggest cyclic flat with respect to the rank that contains $Z_{1}$ and is below a coatom, \ie , if $Z' \in \ZZ(M)$ is such that $Z_{1} \subseteq Z'$ then either $\rho(Z')\geq k-1$ or $\rho(Z') \leq \rho(Z_{m})$. Now, by the same arguments as in the proof of Proposition \ref{prop:ctrapos_double_nulledge}, $M/Z_{m}$ is simple with no isthmuses and has a lattice of cyclic flats of height $3$. Therefore, we have $d \leq d_{M/Z_{m}} \leq 4$. 

Suppose now that there exists $Z_{2} \in \ZZ(M)$ such that $Z_{1} \subset Z_{2}$ and $\rho(Z_{2})=k-2$. We can apply Lemma \ref{lemma:rankdiff2_null_relations} to get a bound on the minimum distance $d$. Let $\Upsilon_{Z_{2}}$ be the set of coatoms containing $Z_{2}$. By Lemma \ref{lemma:upsilon2}, we have $|\Upsilon_{Z_{2}}| \geq 2$ which reduces the possible cases in Lemma \ref{lemma:rankdiff2_null_relations}. We also use the fact that since $\eta(Z_{1}) = \eta(Z^{c})-1$, we have $\eta(Z_{2}) \geq \eta(Z^{c})$.

\begin{enumerate}
\item If $|\Upsilon_{Z_{2}}| = 3$, then we have
\begin{align*}
\eta(E) = 1 + \sum\limits_{Z \in \Upsilon_{Z_{2}}} \eta(Z) - 2 \eta(Z_{2}) \leq 1 + 3 \eta(Z^{d}) - 2 \eta(Z^{c})  & \iff \\
\eta(E) + 1 - \eta(Z^{d}) \leq 2 ( \eta(Z^{d}) - \eta(Z^{c}) + 1) & \iff \\
d \leq 2 ( \eta(Z^{d}) - \eta(Z^{c}) + 1).
\end{align*}
\item If $|\Upsilon_{Z_{2}}|=2$ and $E-(\bigcup\limits_{Z \in \Upsilon_{Z_{2}}} Z) \neq \emptyset$, then we have
\begin{align*}
\eta(E) = 1 + \sum\limits_{Z \in \Upsilon_{Z_{2}}} \eta(Z) - \eta(Z_{2}) \leq 1 + 2 \eta(Z^{d}) - \eta(Z^{c}) &  \iff \\
d \leq 2 + (\eta(Z^{d}) - \eta(Z^{c})).\\
\end{align*}
\item if $|\Upsilon_{Z_{2}}|=2$ and $E=\bigcup\limits_{Z \in \Upsilon_{Z_{2}}} Z$, then we have
\begin{align*}
\eta(E) = \sum\limits_{Z \in \Upsilon_{Z_{2}}} \eta(Z) - \eta(Z_{2}) \leq  2 \eta(Z^{d}) - \eta(Z^{c}) & \iff \\
d \leq 1 + (\eta(Z^{d}) - \eta(Z^{c})).\\
\end{align*}
\end{enumerate}
Hence, the general upper bound for $d$ is the largest of the three bounds obtained above and we have indeed that $d \leq 2(\eta(Z^{d})-\eta(Z^{c}) +1)$.
\end{proof}

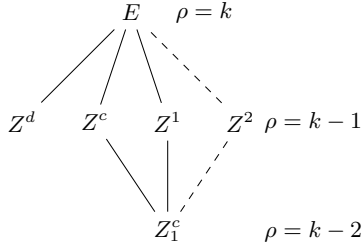
\begin{figure}
\centering
\resizebox{0.4\textwidth}{!}{%
\begin{tikzpicture}
%Starting from top right
\node[shape=rectangle, draw=white] (00) at (1.5,0) {\small $E$};
\node[] (09) at (2.5,0) {\small $\rho=k$};

\node[shape=rectangle,draw=white] (13) at (0,-1.5) {\small $Z^{d}$};
\node[shape=rectangle,draw=white] (10) at (1,-1.5) {\small $Z^{c}$};
\node[shape=rectangle,draw=white] (11) at (2,-1.5) {\small $Z^{1}$};
\node[shape=rectangle,draw=white,dashed] (12) at (3,-1.5) {\small $Z^{2}$};
\node[] (19) at (4,-1.5) {\small $\rho=k-1$};

\node[shape=rectangle,draw=white] (20) at (2,-3) {\small $Z_{1}^{c}$};
\node[] (29) at (4,-3) {\small $\rho=k-2$};
  
\path [-] (10) edge (00);
\path [-] (11) edge (00);
\path [-,dashed] (12) edge (00);
\path [-] (13) edge (00);

\path [-] (20) edge (10);
\path [-] (20) edge (11);
\path [-,dashed] (20) edge (12);

\end{tikzpicture}
}
\caption{Illustration of the hypotheses in Lemma \ref{lemma:dMZc_formulas}. }
\label{fig:lemma_dMZc}
\end{figure}

The strength of this proposition does not really reside in the bound on the minimum distance but in its contrapositive as it gives a sufficient condition for a coatom of $\ZZ(M)$ to be blunt. In order to extend the bluntness property into a recursive structure, we study the minimum distance of $M|Z^{c}$, the matroid restricted to a coatom $Z^{c} \in \ZZ(M)$. We start by a technical lemma which relates $d_{M|Z^{c}}$, $d$, and the nullity of certain coatoms of $\ZZ(M)$. The hypotheses of the next lemma are illustrated in Figure \ref{fig:lemma_dMZc}. 

\begin{lemma}
\label{lemma:dMZc_formulas}
Let $\rmatroid$ be a binary non-degenerate $(n,k,d)$-matroid with $d \geq 3$ and $Z^{d}$ a cyclic flat with maximal nullity. Let $Z^{c} \in \ZZ(M)$ be a blunt coatom with rank $\rho(Z^{c})=k-1$, $Z_{1}^{c} \in \ZZ(M)$ such that $d_{M|Z^{c}}=\eta(Z^{c}) + 1 - \eta(Z_{1}^{c})$, and $\Upsilon_{Z_{1}^{c}}$ the set of coatoms containing $Z_{1}^{c}$. We denote the coatoms in $\Upsilon_{Z_{1}^{c}}$ by $Z^{c}, Z^{1}$ and if it exists, by $Z^{2}$. 

Then, $d_{M|Z^{c}}$ satisfies one of the following.
\begin{enumerate}
\item If $|\Upsilon_{Z_{1}^{c}}|=3$, then $2d_{M|Z^{c}}=d + \eta(Z^{d}) + \eta(Z^{c}) - (\eta(Z^{1}) + \eta(Z^{2}))$.
\item If $|\Upsilon_{Z_{1}^{c}}|=2$ and $E-(Z^{c} \cup Z^{1}) \neq \emptyset$, then $d_{M|Z^{c}} = d - 1 + \eta(Z^{d}) - \eta(Z^{1})$.
\item if $|\Upsilon_{Z_{1}^{c}}|=2$ and $E=Z^{c} \cup Z^{1}$, then we have $d_{M|Z^{c}}= d + \eta(Z^{d}) - \eta(Z^{1})$.
\end{enumerate}
\end{lemma}

\begin{proof}
Since $Z^{c}$ is blunt, we have that $\rho(Z_{1}^{c}) = k-2$ and we can use Lemma \ref{lemma:rankdiff2_null_relations}. We check all possible cases in Lemma \ref{lemma:rankdiff2_null_relations} depending on $\Upsilon_{Z_{1}^{c}}$. Notice that by Lemma \ref{lemma:upsilon2}, we have $|\Upsilon_{Z_{1}}| \geq 2$. 

\begin{enumerate}
\item If $|\Upsilon_{Z_{1}^{c}}| = 3$, then we have
\begin{align*}
\eta(E) = 1 + \eta(Z^{c}) + \sum\limits_{Z \in \Upsilon_{Z_{1}^{c}}-Z^{c}} \eta(Z) - 2 \eta(Z_{1}^{c}) & \iff \\
d_{M|Z^{c}} = \eta(E) - \eta(Z^{d}) + \eta(Z^{d}) - \sum\limits_{Z \in \Upsilon_{Z_{1}^{c}}-Z^{c}} \eta(Z) + \eta(Z^{c}) - \eta(Z^{c}) + \eta(Z_{1}^{c}) & \iff \\
d_{M|Z^{c}} = d - 1 + \eta(Z^{d}) - \sum\limits_{Z \in \Upsilon_{Z_{1}^{c}}-Z^{c}} \eta(Z) + \eta(Z^{c}) - (d_{M|Z^{c}}-1) & \iff \\
2 d_{M|Z^{c}} = d + \eta(Z^{d}) + \eta(Z^{c}) - (\eta(Z^{1}) + \eta(Z^{2})).
\end{align*}
\item If $|\Upsilon_{Z_{1}^{c}}|=2$ and $E-(Z^{c} \cup Z^{1}) \neq \emptyset$, then we have
\begin{align*}
\eta(E) = 1 + \eta(Z^{c}) + \eta(Z^{1}) - \eta(Z_{1}^{c}) & \iff \\
d_{M|Z^{c}} = \eta(E) - \eta(Z^{d}) + \eta(Z^{d}) - \eta(Z^{1}) & \iff \\
d_{M|Z^{c}} = d- 1 + \eta(Z^{d}) - \eta(Z^{1}). \\
\end{align*}
\item if $|\Upsilon_{Z_{1}^{c}}|=2$ and $E=Z^{c} \cup Z^{1}$, then we have
\begin{align*}
\eta(E) = \eta(Z^{c}) + \eta(Z^{1}) - \eta(Z_{1}^{c}) & \iff \\
d_{M|Z^{c}} = d + \eta(Z^{d}) - \eta(Z^{1}). \\
\end{align*}
\end{enumerate}
\end{proof}

From the previous lemma, we can derive a lower bound on $d_{M|Z^{c}}$, which is easier to estimate. 

\begin{proposition}
\label{prop:d_MZd}
Let $\rmatroid$ be a binary non-degenerate $(n,k,d)$-matroid and $Z^{d} \in \ZZ(M)$ a cyclic flat with maximal nullity. If $Z^{c} \in \ZZ(M)$ is a blunt coatom of rank $k-1$, then
\[
d_{M|Z^{c}} \geq \frac{d - ( \eta(Z^{d}) -\eta(Z^{c}) )}{2}.
\]
\end{proposition}

\begin{proof}
Since $Z^{c}$ is a cyclic flat, we have directly that $d_{M|Z^{c}} \geq 2$. Now if $d=2$, we have 
\[
\frac{2-( \eta(Z^{d}) -\eta(Z^{c}))}{2} \leq 1 \leq d_{M|Z^{c}}.
\]
If $d \geq 3$, this proposition is a direct consequence of the previous Lemma \ref{lemma:dMZc_formulas}. Namely, for all $Z \in \ZZ(M)-E$ we have $\eta(Z) \leq \eta(Z^{d})$. By replacing the unknown nullities in Lemma \ref{lemma:dMZc_formulas} with $\eta(Z^{d})$, we get three lower bounds on $d_{M|Z^{d}}$. Therefore, the general bound is the smallest lower bound, which is when $|\Upsilon_{Z_{1}^{c}}|=3$ and $d_{M|Z^{c}} \geq \frac{d-( \eta(Z^{d}) -\eta(Z^{c}))}{2}$.
\end{proof}

The contrapositives of Propositions \ref{prop:ctrapos_double_nulledge} and \ref{prop:bound_null_edge_d} reveals how the minimum distance forces many coatoms $Z^{c} \in \ZZ(M)$ to be blunt. Now, given the lower bound on the minimum distance $d_{M|Z^{c}}$ provided by Proposition \ref{prop:d_MZd}, we can apply again Propositions \ref{prop:ctrapos_double_nulledge} and \ref{prop:bound_null_edge_d} to the restricted matroid $M|Z^{c}$ leading to more blunt cyclic flats. By repeating this process, we obtain decreasing chains of blunt cyclic flats with upper bounded nullity in $\ZZ(M)$. The next example illustrates the strength of Propositions \ref{prop:bound_null_edge_d} and \ref{prop:d_MZd} for the study of the lattice of cyclic flats together with specific techniques on the relation between the nullity and the minimum distance. 

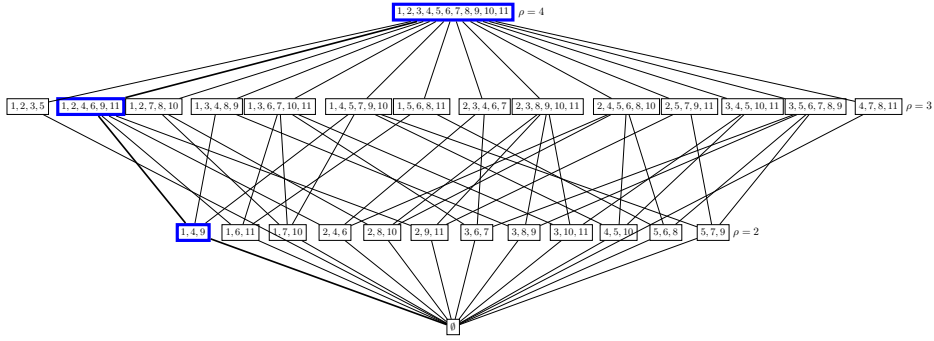
\begin{figure}
\centering
\resizebox{0.98\textwidth}{!}{%
\begin{tikzpicture}

\node[shape=rectangle, draw=black] (00) at (0,0) {\small$\emptyset$};

\node[shape=rectangle,draw=blue, line width = 1mm] (10) at (-8.25,3) {\small\(1,4,9\)};
\node[shape=rectangle,draw=black] (11) at (-6.75,3) {\small\(1,6,11\)};
\node[shape=rectangle,draw=black] (12) at (-5.25,3) {\small\(1,7,10\)};
\node[shape=rectangle,draw=black] (13) at (-3.75,3) {\small\(2,4,6\)};
\node[shape=rectangle,draw=black] (14) at (-2.25,3) {\small\(2,8,10\)};
\node[shape=rectangle,draw=black] (15) at (-0.75,3) {\small\(2,9,11\)};
\node[shape=rectangle,draw=black] (16) at (0.75,3) {\small\(3,6,7\)};
\node[shape=rectangle,draw=black] (17) at (2.25,3) {\small\(3,8,9\)};
\node[shape=rectangle,draw=black] (18) at (3.75,3) {\small\(3,10,11\)};
\node[shape=rectangle,draw=black] (19) at (5.25,3) {\small\(4,5,10\)};
\node[shape=rectangle,draw=black] (110) at (6.75,3) {\small\(5,6,8\)};
\node[shape=rectangle,draw=black, label=right:{$\rho = 2$}] (111) at (8.25,3) {\small\(5,7,9\)};
 
\node[shape=rectangle,draw=black] (20) at (-13.5,7) {\small\(1,2,3,5\)};
\node[shape=rectangle,draw=blue, line width = 1mm] (21) at (-11.5,7) {\small\(1,2,4,6,9,11\)};
\node[shape=rectangle,draw=black] (22) at (-9.5,7) {\small\(1,2,7,8,10\)};
\node[shape=rectangle,draw=black] (23) at (-7.5,7) {\small\(1,3,4,8,9\)};
\node[shape=rectangle,draw=black] (24) at (-5.5,7) {\small\(1,3,6,7,10,11\)};
\node[shape=rectangle,draw=black] (25) at (-3,7) {\small\(1,4,5,7,9,10\)};
\node[shape=rectangle,draw=black] (26) at (-1,7) {\small\(1,5,6,8,11\)};
\node[shape=rectangle,draw=black] (27) at (1,7) {\small\(2,3,4,6,7\)};
\node[shape=rectangle,draw=black] (28) at (3,7) {\small\(2,3,8,9,10,11\)};
\node[shape=rectangle,draw=black] (29) at (5.5,7) {\small\(2,4,5,6,8,10\)};
\node[shape=rectangle,draw=black] (210) at (7.5,7) {\small\(2,5,7,9,11\)};
\node[shape=rectangle,draw=black] (211) at (9.5,7) {\small\(3,4,5,10,11\)};
\node[shape=rectangle,draw=black] (212) at (11.5,7) {\small\(3,5,6,7,8,9\)};
\node[shape=rectangle,draw=black, label=right:{$\rho = 3$}] (213) at (13.5,7) {\small\(4,7,8,11\)};
 
 \node[shape=rectangle,draw=blue, line width=1mm, label=right:{$\rho = 4$}] (30) at (0,10) {\small\(1,2,3,4,5,6,7,8,9,10,11\)};

\path [-, draw=blue, line width=1.5pt] (00) edge (10);
\path [-] (00) edge (11);
\path [- ] (00) edge (12);
\path [- ] (00) edge (13);
\path [- ] (00) edge (14);
\path [- ] (00) edge (15);
\path [- ] (00) edge (16);
\path [- ] (00) edge (17);
\path [- ] (00) edge (18);
\path [- ] (00) edge (19);
\path [- ] (00) edge (110);
\path [- ] (00) edge (111);
 
\path [- ] (00) edge (20);
\path [- ] (00) edge (213);
 \path [-,draw = blue, line width = 1.5pt] (10) edge (21);
\path [-] (10) edge (23);
\path [-] (10) edge (25);
 \path [-] (11) edge (21);
\path [-] (11) edge (24);
\path [-] (11) edge (26);
 \path [-] (12) edge (22);
\path [-] (12) edge (24);
\path [-] (12) edge (25);
 \path [-] (13) edge (21);
\path [-] (13) edge (27);
\path [-] (13) edge (29);
 \path [-] (14) edge (22);
\path [-] (14) edge (28);
\path [-] (14) edge (29);
 \path [-] (15) edge (21);
\path [-] (15) edge (28);
\path [-] (15) edge (210);
 \path [-] (16) edge (24);
\path [-] (16) edge (27);
\path [-] (16) edge (212);
 \path [-] (17) edge (23);
\path [-] (17) edge (28);
\path [-] (17) edge (212);
 \path [-] (18) edge (24);
\path [-] (18) edge (28);
\path [-] (18) edge (211);
 \path [-] (19) edge (25);
\path [-] (19) edge (29);
\path [-] (19) edge (211);
 \path [-] (110) edge (26);
\path [-] (110) edge (29);
\path [-] (110) edge (212);
 \path [-] (111) edge (25);
\path [-] (111) edge (210);
\path [-] (111) edge (212);

 \path [-] (20) edge (30);
\path [-,draw = blue, line width = 1.5pt] (21) edge (30);
\path [-] (22) edge (30);
\path [-] (23) edge (30);
\path [-] (24) edge (30);
\path [-] (25) edge (30);
\path [-] (26) edge (30);
\path [-] (27) edge (30);
\path [-] (28) edge (30);
\path [-] (29) edge (30);
\path [-] (210) edge (30);
\path [-] (211) edge (30);
\path [-] (212) edge (30);
\path [-] (213) edge (30);

\end{tikzpicture}
}
\caption{Lattice of cyclic flats of the binary $(11,4,5)$-matroid from Example \ref{ex:110405}. }
\label{fig:1104_blue}
\end{figure}

\setcounter{MaxMatrixCols}{20}
\begin{example}
\label{ex:110405}
Let $\rmatroid$ be a binary simple $(11,4,5)$-matroid. Since the minimum distance is equal to $5$, we know that $E$ is blunt and all coatoms have rank $k-1$. Moreover, there is a coatom $Z^{d} \in \ZZ(M)$ with size $6$ and rank $3$. By Proposition \ref{prop:ctrapos_double_nulledge}, $Z^{d}$ is blunt and $\eta(Z^{d})=3 > 1$. Proposition \ref{prop:d_MZd} implies that $d_{M|Z^{d}} \geq \left\lceil \frac{d}{2} \right\rceil = 3$. Since the maximal minimum distance for a $(6,3,d_{M|Z^{d}})$ matroid is $3$, we have directly that $d_{M|Z^{d}}=3$. Now by Theorem \ref{thm:classif_height3}, $M|Z^{d}$ is isomorphic to the $(6,3,3)$-matroid studied in Section \ref{sec:height3}. Therefore, there are exactly $4$ cyclic flats contained in $Z^{d}$ with size $3$ and rank $2$. 

Let $Z_{1}^{d} \in \ZZ(M)$ be such that $Z_{1}^{d} \lessdot Z^{d}$. We can apply Lemma \ref{lemma:dMZc_formulas} to obtain the nullity of the other coatoms containing $Z_{1}^{d}$. Indeed, since $d_{M|Z^{d}}=3<d-1=4$, we know by Lemma \ref{lemma:dMZc_formulas} that $\Upsilon_{Z_{1}^{d}}$, the set of coatoms containing $Z_{1}^{d}$, has size $|\Upsilon_{Z_{1}^{d}}|=3$. Let $Z^{1}$ and $Z^{2}$ be the two other coatoms in $\Upsilon_{Z_{1}^{d}}$. The first formula in Lemma \ref{lemma:dMZc_formulas} simplifies as $\eta(Z^{1}) + \eta(Z^{2}) = 5$. Now the maximal nullity of a coatom is $\eta(Z^{d})=3$. Thus we have $\eta(Z^{1})=3$ and $\eta(Z^{2})=2$. Hence, there exist coatoms with parameters $(5,3,d')$. By Proposition \ref{prop:bound_null_edge_d}, the minimum distance $d_{M|Z^{2}}$ is at least $2$ and since it contains already a cyclic flat with nullity 1, we have that $d_{M|Z^{2}}=2$. 

In summary, $M$ contains at least $5$ different cyclic flats $(6,3,3)$, $4$ cyclic flats $(5,3,2)$ and $8$ atoms $(2,1,2)$. It is, in fact, possible to obtain the remaining cyclic flats by using some arguments about the intersection between two coatoms but this is rather long and mostly specific to this particular example. We can now double check our results by finding a particular generator matrix for $M$ and displaying the lattice of cyclic flats. Let $G$ be the following matrix : 
\[
G=
\begin{pmatrix}
1&0&0&0&1&0&0&1&1&1&1\\
0&1&0&0&1&1&1&0&0&1&1\\
0&0&1&0&1&0&1&1&0&1&0\\
0&0&0&1&0&1&1&1&1&1&1\\
\end{pmatrix}.
\]
The matroid $M(G)$ is indeed a simple $(11,4,5)$-matroid and its lattice of cyclic flats is displayed in Figure \ref{fig:1104_blue}.
\end{example}

\subsection{Residual Codes and the Griesmer Bound}

The final part of this section is dedicated to reformulating two notions in coding theory known as residual codes and the Griesmer bound for binary matroids. For more information about these two notions, we refer the reader to \cite[Section~2.7]{huffman2010}. The main result is an extension of Proposition \ref{prop:d_MZd} to arbitrary binary matroids which is the exact correspondent of the existence of residual codes for binary linear codes. As a direct consequence, we obtain the Griesmer bound for binary matroids. 

\begin{theorem}
\label{thm:resid_code}
If $\rmatroid$ is a binary $(n,k,d)$-matroid, then there exists $A\subset E$ such that $M|A$ is a binary $(n-d,k-1,d')$-matroid with $d'\geq \frac{d}{2}$.
\end{theorem}

\begin{corollary}
\label{cor:griesmer_bound}
For a binary $(n,k,d)$-matroid, we have
\[
n \geq \sum\limits_{i=0}^{k-1}\left\lceil \frac{d}{2^{i}} \right\rceil.
\]
\end{corollary}

We start by the proof of Theorem \ref{thm:resid_code}. 

\begin{proof}
Let $\rmatroid$ be a binary $(n,k,d)$-matroid. We separate the proof into three cases in which we give an explicit construction of the set $A$ with the required parameters. Notice that the covering relations of the lattice $\ZZ(M)$ are not affected by the existence of loops nor is the minimum distance since the nullity of all cyclic flats increases evenly by the number of loops. Thus, without loss of generality, we assume that $M$ contains no loops. Notice also that if $M$ contains no isthmuses then the restriction to a cyclic flat will not create any isthmuses. 

\begin{enumerate}
\item Suppose $M$ contains no isthmuses and $d \geq 3$. Let $Z^{d} \in \ZZ(M)$ with $|Z^{d}| = n-d$. 
\begin{itemize}
\item If $Z^{d}$ is blunt, then Proposition \ref{prop:d_MZd} guarantees that $d_{M|Z^{d}} \geq \frac{d}{2}$. Thus, we can choose $A=Z^{d}$.  
\item If $Z^{d}$ is not blunt, then this implies that $d_{M|Z^{d}} =2$. By Proposition \ref{prop:ctrapos_double_nulledge}, we have $d \leq 4$. So indeed $d_{M|Z^{d}} \geq \frac{d}{2}$ and we can choose $A=Z^{d}$. 
\end{itemize}
\item Suppose $M$ contains no isthmuses and $d=2$. 
\begin{itemize}
\item If there exists $Z^{d} \in \ZZ(M)$ with maximal nullity and $\rho(Z^{d}) = k-1$, then $|Z^{d}| = n-d$ and $d_{M|Z^{d}} \geq 2$. So indeed $d_{M|Z^{d}} \geq \frac{d}{2}=1$ and we can choose $A=Z^{d}$.
\item Assume there is no such $Z^{d}$. Let $Z$ be such that $Z \lessdot E$ and $\eta(Z)=\eta(E) + 1 -d = \eta(E)-1$. Now $M/Z$ is isomorphic to a uniform matroid $U_{m}^{m+1}$ with $m=k-\rho(Z)$. This implies that if $B \subset E-Z$ with $|B|=m-1$, then $\rho(Z\cup B) = \rho(Z) + |B| = k-1$. We also have that $|Z\cup B |= |Z| + |B| = \rho(Z) + \eta(Z) + m-1 = k - 1 + \eta(E) - 1 = n-2 = n-d$. Finally, the minimum distance is $d_{M|Z\cup B} = 1 = \frac{d}{2}$. Thus, we can choose $A=Z\cup B$. 
\end{itemize}
\item Finally, suppose $M$ contains some isthmuses. This implies that the minimum distance $d$ is equal to 1. Let $H$ be a hyperplane of $M$. $H$ has parameters $\rho(H)=k-1$, $|H|=n-1$, and $d_{M|H} \geq 1 \geq \frac{d}{2}$. Thus, we can choose $A=H$. 
\end{enumerate}
\end{proof}

We give a proof of Corollary \ref{cor:griesmer_bound} for completeness. This proof follows a standard proof of the Griesmer bound by using residual codes as in \cite[Theorem~2.7.4]{huffman2010}.  

\begin{proof}
Notice first that if $d_{M|Z} \geq \frac{d}{2}$ then $d_{M|Z} \geq \left\lceil \frac{d}{2} \right\rceil$. We will now prove the statement by induction on $k$. 

If $k=1$, then the conclusion is trivial since it says that $n \geq d$. Let $k>1$ and assume that the statement is true for any binary matroid with rank $k-1$. By Theorem \ref{thm:resid_code}, there exists a subset $A \subset E$ such that $M|A$ has size $n-d$, dimension $k-1$ and minimum distance $d_{M|A} \geq \left\lceil \frac{d}{2} \right\rceil$. By applying the induction hypothesis on $M|A$, we have
\[
n-d \geq \sum\limits_{i=0}^{k-2}\left\lceil \frac{d_{M|A}}{2^{i}} \right\rceil \geq \sum\limits_{i=0}^{k-2} \left\lceil \frac{\left\lceil \frac{d}{2} \right\rceil}{2^{i}} \right\rceil = \sum\limits_{i=0}^{k-2}\left\lceil \frac{d}{2^{i+1}} \right\rceil .
\]
Now, we add $d$ to both sides to get
\[
n \geq \sum\limits_{i=0}^{k-1}\left\lceil \frac{d}{2^{i}} \right\rceil.
\]
\end{proof}

By combining the two previous proofs, we can understand the Griesmer bound as an evaluation of the parameters of a chain contained almost entirely in $\ZZ(M)$. Indeed, every subset $A \subset E$ that we constructed in the proof of Theorem \ref{thm:resid_code} is a flat if not directly a cyclic flat. Furthermore, since $\ZZ(M|F)$ is a sub-lattice of $\ZZ(M)$ when $F$ is a flat, performing the recursive steps of choosing residual codes can be viewed as taking a decreasing chain in the lattice of cyclic flats completed by the lattice of flats for every encounter of a rank edge. Finally, as illustrated in the next example, the construction of such a chain can be directly extracted from the proof of Theorem \ref{thm:resid_code}.

\begin{example}
Let $M$ be the binary $(11,4,5)$-matroid given in Example \ref{ex:110405}. $M$ achieves the Griesmer bound since we have $11 = \sum_{i=0}^{3}\left\lceil \frac{5}{2^{i}} \right\rceil = 5 + 3 + 2 + 1$. Now we construct a decreasing chain $E\gtrdot Z^{d} \gtrdot Z_{1}^{d} \gtrdot \emptyset$ in $\ZZ(M)$ by taking at every step a cyclic flat with maximal nullity contained in the previous one. By labelling the columns of the generator matrix $G$ from $1$ to $11$, one such chain is given by $[11] \gtrdot \{1,2,4,6,9,11\} \gtrdot \{1,4,9\} \gtrdot \emptyset$ and is displayed in blue in Figure \ref{fig:1104_blue}. 

In Example \ref{ex:110405}, we saw that $M|Z^{d}$ is a $(6,3,3)$-matroid and $M|Z_{1}^{d}$ is a $(3,2,2)$-matroid. Since $\emptyset \lessdot Z_{1}^{d}$ is a rank edge, we complete the chain by adding the flat $\{ e \} \in \FF(M)$ with $e \in Z_{1}^{d}$. Hence we have $n=d+d_{M|Z^{d}} + d_{M|Z_{1}^{d}} + d_{M|\{e\}} = 5 + 3 + 2 + 1 = \sum\limits_{i=0}^{3}\left\lceil \frac{5}{2^{i}} \right\rceil$. 
\end{example}

The previous proofs give a new understanding of the Griesmer bound and they are, in fact, deeply connected with the standard proofs in coding theory. The existence of residual codes is usually proven by using a codeword with desired weight and puncturing on its support \cite{huffman2010}. Therefore, to show the link between the different proofs, we demonstrate the relation between codewords of a binary linear code and coatoms of the lattice of cyclic flats in the associated matroid. 

\begin{lemma}
\label{lemma:codewords_to_coatoms}
Let $\CC$ be a binary non-degenerate $[n,k,d]$ linear code with $d\geq 3$ and let $c$ be a codeword of $\CC$ with weight $\wt(c)<2d-2$. Then, $Z_{c}=E-\supp(c)$ is a coatom of $\ZZ(M)$. 
\end{lemma}

\begin{proof}
Let $S=\supp(c)$. We have that $|Z_{c}|=|E|-|S| = n-\wt(c)$. Assume for a contradiction that $\rho(Z_{c})<k-1$. Then there exists $c'$ a codeword of $\CC$ different from $c$ such that $c'_{i}=0$ for all $i \in Z_{c}$. Now let $\alpha \in \F_{2}$ such that at least $\wt(c)/2$ coordinates of $c'_{|S}$ equal $\alpha$. Then, we have
\[
d \leq \wt(c'-\alpha c) \leq \wt(c) - \frac{\wt(c)}{2} = \frac{\wt(c)}{2}
\]
which contradicts the hypothesis on $\wt(c)$. Thus, the rank of $Z_{c}$ is $k-1$. 

$Z_{c}$ is a flat, since otherwise, $c$ does not have weight $\wt(c)$. It remains to prove that $Z_{c}$ is cyclic. Assume for a contradiction that there exists $e \in Z_{c}$ such that $\rho(Z_{c}-\{e \} ) < \rho(Z_{c})$. Since there is an isthmus in $M|Z_{c}$, this implies that $d_{M|Z_{c}}=1$ and there is a codeword $\hat{c}$ such that $\wt(\hat{c}_{|Z_{c}})=1$. Let $\beta \in \F_{2}$ be such that at least $\wt(c)/2$ coordinates of $\hat{c}_{|S}$ equal $\beta$. Then, we have
\[
d \leq \wt(\hat{c} - \beta c) \leq \wt(c) + 1 - \frac{\wt(c)}{2} = \frac{\wt(c) + 2}{2} < \frac{2d - 2 + 2}{2} = d
\]
which is a contradiction. Hence, $Z_{c}$ is a cyclic flat of rank $k-1$ and thus a coatom of $\ZZ(M)$. 
\end{proof}

\begin{lemma}
Let $M$ be a binary non-degenerate $(n,k,d)$-matroid with $d\geq 3$. Let $Z_{c}$ be a coatom of $\ZZ(M)$ with $|Z_{c}|>n-2d+2$. Then, there exists a codeword $c$ in $\CC_{M}$, the linear code associated to $M$, such that $\supp(c) = E - Z_{c} $. 
\end{lemma}

\begin{proof}
Let $G_{M}$ be a generator matrix of $M$. Since $d \geq 3$, we have $\rho(Z^{d})=k-1$ and in particular, $G_{M|Z^{d}}$ the submatrix of $G_{M}$ restricted to the columns indexed by $Z^{d}$, has rank $k-1$. Since $G_{M}$ has $k$ rows, one of the rows in $G_{M|Z^{d}}$ is dependent of the others. This implies that there is a codeword $c$ in $\CC_{M}$ such that $c_{j}=0$ for all $j \in Z^{d}$. Since $|Z^{d}|>n-2d+2$, we have $\wt(c)< 2d+2$ and $\supp(c)\subseteq E-Z_{c}$. By Lemma \ref{lemma:codewords_to_coatoms}, $E-\supp(c)$ is a coatom of $\ZZ(M)$ and $Z_{c} \subseteq E-\supp(c)$. Hence, $Z_{c}=E-\supp(c)$ since $Z_{c}$ is already a coatom of $\ZZ(M)$. 
\end{proof}

By combining the two previous lemmas, we get the following result. 

\begin{proposition}
Let $\CC$ be a non-degenerate binary $[n,k,d]$ linear code and $M_{C}$ the associated matroid. Then, there is a bijective map between the codewords of weight less than $2d-2$ and the coatoms of $\ZZ(M_{C})$ of size less than $n-2d+2$. 
\end{proposition}

Thus, this result shows the relation between the small weight codewords of a binary linear code and the cyclic flats with a small size of the matroid associated to the code. 

 % Recursive structure and Griesmer Bound
%Input "conclusions"
%Contains : What do you think!?

\section{Conclusion}

In this paper, we presented the first steps towards the characterization of the lattice of cyclic flats of representable matroids over $\F_{q}$. In the first part of the paper, we derived two natural maps from $\ZZ(M)$ to the lattice of cyclic flats of a minor $M|Y/X$. Then, we showed how to reconstruct the lattice of flats from $\ZZ(M)$ and we computed the largest $n$ for which the uniform matroid $U_{n}^{2}$ is a minor of $M$, from the lattice of cyclic flats. In the second part, we focused on binary matroids and the structure of their lattice of cyclic flats. We proved that the lattice of cyclic flats of a binary simple matroid with no isthmuses is atomic. Furthermore, we classified the binary matroids with lattice of cyclic flats of height $3$. Finally, we defined the class of blunt cyclic flats for binary matroids and demonstrated the relation between blunt cyclic flats and the minimum distance of a matroid. As a consequence of this relation, we reproved the Griesmer bound for binary codes. % Conclusion

\section{Acknowledgements}
The work of M. Grezet, C. Hollanti, and T.  Westerb\"ack was supported in part by the Academy of Finland [grant numbers 276031, 282938, and 303819] and by the Technical University of Munich – Institute for Advanced Study, funded by the German Excellence Initiative and the EU 7th Framework Programme [grant number 291763], via a Hans Fischer Fellowship.
R. Freij-Hollanti was supported by the German Research Foundation (Deutsche Forschungsgemeinschaft, DFG) [grant number WA3907/1-1].

\bibliographystyle{elsarticle-num}
\bibliography{references}

\begin{thebibliography}{10}
\expandafter\ifx\csname url\endcsname\relax
  \def\url#1{\texttt{#1}}\fi
\expandafter\ifx\csname urlprefix\endcsname\relax\def\urlprefix{URL }\fi
\expandafter\ifx\csname href\endcsname\relax
  \def\href#1#2{#2} \def\path#1{#1}\fi

\bibitem{birkhoff1995}
G.~Birkhoff, Lattice Theory, 3rd Edition, Vol.~25 of Colloquium Publications,
  American Mathematical Society, 1967.

\bibitem{sims80}
J.~A. Sims, Some problems in matroid theory, Ph.D. thesis, University of Oxford
  (1980).

\bibitem{bonin08}
J.~E. Bonin, A.~de~Mier, The lattice of cyclic flats of a matroid, Annals of
  Combinatorics 12~(2) (2008) 155--170.

\bibitem{prideaux16}
K.~Prideaux, Matroids, cyclic flats, and polyhedra, Master's thesis, Victoria
  University of Wellington (2016).

\bibitem{eberhardt14}
J.~N. Eberhardt, Computing the {T}utte polynomial of a matroid from its lattice
  of cyclic flats, The Electronic Journal of Combinatorics 21~(3) (2014) 3--47.

\bibitem{westerback15}
T.~Westerb{\"a}ck, R.~Freij-Hollanti, T.~Ernvall, C.~Hollanti, On the
  combinatorics of locally repairable codes via matroid theory, IEEE
  Transactions on Information Theory 62~(10) (2016) 5296--5315.

\bibitem{rota71}
G.-C. Rota, Combinatorial theory, old and new, in: Actes du Congr\`{e}s
  International des Math\'{e}maticiens (Nice, 1970), 1971, pp. 229--233.

\bibitem{whittle14}
J.~Geelen, B.~Gerards, G.~Whittle, Solving {R}ota's conjecture, Notices of the
  American Mathematical Society 61~(7) (2014) 736--743.

\bibitem{oxley11}
J.~Oxley, Matroid Theory, 2nd Edition, Vol.~21 of Oxford Graduate Texts in
  Mathematics, Oxford University Press, 2011.

\bibitem{papailiopoulos12}
D.~S. Papailiopoulos, A.~G. Dimakis, Locally repairable codes, in: Proceedings
  of the IEEE International Symposium on Information Theory, 2012, pp.
  2771--2775.

\bibitem{silberstein18}
N.~Silberstein, A.~Zeh, Anticode-based locally repairable codes with high
  availability, Designs, Codes and Cryptography 86~(2) (2018) 419--445.

\bibitem{huang15}
P.~Huang, E.~Yaakobi, H.~Uchikawa, P.~H. Siegel, Cyclic linear binary locally
  repairable codes, in: Proceedings of the IEEE information theory workshop,
  IEEE, 2015, pp. 1--5.

\bibitem{huffman2010}
W.~C. Huffman, V.~Pless, Fundamentals of {E}rror-{C}orrecting {C}odes,
  Cambridge University Press, 2010.

\bibitem{grezet17}
M.~Grezet, R.~Freij-Hollanti, T.~Westerb\"ack, C.~Hollanti, On binary matroid
  minors and applications to data storage over small fields, in: International
  Castle Meeting on Coding Theory and Applications, 2017, pp. 139--153.

\bibitem{grezet18}
M.~Grezet, R.~Freij-Hollanti, T.~Westerb{\"a}ck, O.~Olmez, C.~Hollanti, Bounds
  on binary locally repairable codes tolerating multiple erasures, in: The
  International Zurich Seminar on Information and Communication (IZS).
  Proceedings, Zurich, Switzerland, 2018, pp. 103--107.

\bibitem{vamos78}
P.~V{\'a}mos, The missing axiom of matroid theory is lost forever, Journal of
  the London Mathematical Society 2~(3) (1978) 403--408.

\bibitem{mayhew18}
D.~Mayhew, M.~Newman, G.~Whittle, Yes, the ``missing axiom" of matroid theory
  is lost forever, Transactions of the American Mathematical Society 370~(8)
  (2018) 5907--5929.

\bibitem{tutte58}
W.~T. Tutte, A homotopy theorem for matroids, {I}, {II}, Transactions of the
  American Mathematical Society 88~(1) (1958) 144--174.

\bibitem{Segre55}
B.~Segre, Curve razionali normali e $k$-archi negli spazi finiti, Annali di
  Matematica Pura ed Applicata 39~(1) (1955) 357--379.

\bibitem{BallMDS}
S.~Ball, On sets of vectors of a finite vector space in which every subset of
  basis size is a basis, Journal of the European Mathematical Society 14~(3)
  (2012) 733--748.

\bibitem{crapo70}
H.~H. Crapo, G.-C. Rota, On the Foundations of Combinatorial Theory:
  Combinatorial Geometries, MIT Press, 1970.

\bibitem{stanley2011}
R.~P. Stanley, Enumerative Combinatorics, 2nd Edition, Vol.~1, Cambridge
  University Press, 2011.

\bibitem{shoda16}
K.~Shoda, Large families of matroids with the same {T}utte polynomial, Ph.D.
  thesis, George Washington University (2012).

\end{thebibliography}

\end{document}